\documentclass[review]{elsarticle}

\usepackage{lineno,hyperref}
\modulolinenumbers[1]

\journal{Transportation Research Part B: Methodological}

\biboptions{authoryear}

\usepackage[margin=1in]{geometry}
\usepackage{color,soul}
\usepackage{enumitem}
\usepackage{booktabs}
\usepackage{array}
\newcolumntype{P}[1]{>{\centering\arraybackslash}p{#1}}
\newcolumntype{L}[1]{>{\raggedright\arraybackslash}p{#1}}
\usepackage{rotating}
\usepackage{multirow}
\usepackage{subfigure}
\usepackage{amsmath}
\usepackage{amsfonts}
\usepackage{gentium}
\usepackage{enumitem}
\usepackage{lmodern}
\usepackage{hyperref}
\hypersetup{
     colorlinks   = true,
     citecolor    = gray
}

\begin{document}

\begin{frontmatter}

\title{Advancements in Continuum Approximation Models for Logistics and Transportation Systems: 1996 - 2016}
%\tnotetext[mytitlenote]{Fully documented templates are available in the elsarticle package on \href{http://www.ctan.org/tex-archive/macros/latex/contrib/elsarticle}{CTAN}.}

%% Group authors per affiliation:
%\author{Elsevier\fnref{myfootnote}}
%\address{Radarweg 29, Amsterdam}
%\fntext[myfootnote]{Since 1880.}

%% or include affiliations in footnotes:
\author[northwestern]{Sina Ansari}
\author[northwestern]{Mehmet Ba\c{s}dere\corref{corAuth}}
\cortext[corAuth]{Corresponding author}
\ead{mbasdere@u.northwestern.edu }
\author[usf]{Xiaopeng Li}
\author[uiuc]{Yanfeng Ouyang}
\author[northwestern]{Karen Smilowitz}

\address[northwestern]{Department of Industrial Engineering and Management Sciences, Northwestern University, Evanston, IL 60208}
\address[usf]{Department of Civil and Environmental Engineering, University of South Florida, Tampa, FL 33620}
\address[uiuc]{Department of Civil and Environmental Engineering, University of Illinois at Urbana-Champaign, Urbana, IL 61801}

\begin{abstract}
Continuum Approximation (CA) is an efficient and parsimonious technique for modeling complex logistics problems. In this paper, we review recent studies that develop CA models for transportation, distribution and logistics problems with the aim of synthesizing recent advancements and identifying current research gaps. This survey focuses on important principles and key results from CA models. In particular, we consider how these studies fill the gaps identified by the most recent literature reviews in this field. We observe that CA models are used in a wider range of applications, especially in the areas of facility location and integrated supply chain management. Most studies use CA as an alternative to exact solution approaches; however, CA can also be used in combination with exact approaches. We also conclude with promising areas of future work.
\end{abstract}


\end{frontmatter}

%\linenumbers

\section{Introduction}

%[Basic problem types]

\citet{simchi2008designing} define supply chain and logistics management as \textquotedblleft the set of approaches utilized to efficiently integrate suppliers, manufacturers, warehouses, and stores, so that merchandise is produced and distributed at the right quantities, to the right locations, and at the right time, in order to minimize system wide costs while satisfying service level requirements''. Numerous studies have been conducted to characterize, analyze and optimize planning, design and operations of logistics and transportation systems. Typical examples of such problems include those related to facility location planning and vehicle routing. Traditional approaches tended to characterize these problems in a discrete setting, e.g., with a fixed set of candidate facility locations, discrete time periods, and discrete customer demand points, so that these problems can be solved by well-developed integer mathematical programming techniques. For example, \citet{daskin95} and \citet{drezner95} systematically introduced a range of classic discrete facility location models including covering problems \citep{Christofides75,Church74}, center and median problems \citep{Hakimi64} and fixed-charge location problems \citep{Cornuejols77, Mirzain85}. Later, a series of new discrete models have been proposed to address location problems with stochastic demand \citep{daskin82, daskin83, batta89, Dasci2005} and unreliable facility services\citep{Church74, snyder05, qi2007supply, Berman07, qi2009continuous, ouyang10cui, lim2010facility, chen2011joint, li2011reliable, li2012reliable, yun2015reliability}. Numerous discrete models have also been developed to address vehicle routing issues at the operational level in both deterministic and stochastic environments. See \citet{baldacci2007recent, cordeau2007vehicle, laporte2009fifty, toth2002vehicle, gendreau1996stochastic} for some reviews.
%Challenges to discrete models

Although discrete models, especially with the help of modern computation power, can sometimes yield exact solutions to large-scale logistics problems, they generally have a relatively complex formulation structure that may hinder our understanding of problem properties and managerial insights. Often, the problems belong to the class of NP-hard problems, and hence solving large-scale instances would require enormous computational efforts which likely increase exponentially with the problem instance size. Hence, it is often not practical to solve large-scale logistics problems to optimality. Further, there is often uncertainty in the corresponding data and the lack of precision leads to inaccuracies in the optimal solution \citep{daganzo1987increasing}. These drawbacks are particularly prominent if one attempts to make decisions (e.g., those on location, inventory and routing) in stochastic, time-varying, competitive and coupled environments. For example, stochasticity could arise from both the demand side (e.g., random customers) and the supply side (e.g., service disruptions) and imposes a large number of induced realization scenarios. System operation characteristics, such as link travel time and resource availability, can be time-dependent due to exogenous (e.g., weather condition) or endogenous factors (e.g., congestion). Competition among service providers and/or customers may require equilibrium perspectives to be blended through a hierarchical modeling structure, such as a mathematical program with equilibrium constraints (MPEC) or an equilibrium problem with equilibrium constraints (EPEC). Emerging vehicle technologies (e.g., electric vehicles and autonomous cars) and transportation modes (e.g., car-sharing and ride-sourcing) pose new constraints to daily operations of vehicle fleets (e.g., electric vehicle charging) while creating new mobility paradigms bridging traditional public and private transportation services. 

%Intro to CA

The concept of CA as a complement to discrete models has been shown suitable for addressing these above-mentioned challenges in various contexts. The CA approach was first proposed by \citet{newell71} and \citet{newell73} and has been widely applied to various logistics problems including facility location, inventory management and vehicle routing. CA models feature continuous representations of input data and decision variables as density functions over time and space, and the key idea is to approximate the objective into a functional (e.g., integration) of localized functions that can be optimized by relatively simple analytical operations. Each localized function approximates the cost structure of a local neighborhood with nearly homogeneous settings. Such homogeneous approximation enables mapping otherwise high-dimensional decision variables into a low-dimensional space, allowing the optimal design for this neighborhood to be obtained with simple calculus, even when spatial stochasticity, temporal dynamics and other operational complexities are present. The results from such models often bear closed-form analytical structures that help reveal managerial insights. Compared with their discrete counterparts, CA methods generally incur less computational burden, require less accurate input data, and, more importantly, can conveniently reveal managerial insights, especially for large-scale practical problems. These appealing features have motivated researchers to explore simple solutions for various complex problems arising in the logistics and transportation fields in the past few decades. 

%[Re-organize this paragraph by issues.. summarize work prior to 1996 -- XP]

CA has been applied to three basic logistics problem classes: location, routing and inventory management. In earlier applications, CA was used to determine facility locations and corresponding assignments of customers to these facilities in a continuous space \citep{newell73,daganzo86}. The key to a CA location problem is to balance the tradeoff between long-term transportation cost and one-time facility investment, which is usually formulated as analytical functions of local facility density (or its inverse, a facility's service area). Further, CA is used to formulate routing problems that determine the most economic routes for vehicles to deliver or pickup commodities or people across a continuous space. The scope of routing problems includes single vehicle delivery (also known as the traveling salesman problem) \citep{daganzo84a}, multi-vehicle based distribution \citep{newell1986design}, and multi-echelon distribution with intermediate consolidation and transshipment facilities \citep{daganzo1988comparison}. A fundamental problem in CA-based routing is to format or partition the space into certain geometries suitable for constructing near-optimum vehicle routes with simple heuristics. An inventory management problem investigates the trade-off between the inventory size and the corresponding transportation cost at a supply chain facility \citep{blumenfeld1991synchronizing}. With homogeneous approximation in local spatiotemporal neighborhoods, the basic system cost in an inventory management problem can be often formulated into an economic-order-quantity (EOQ) function that has a simple analytical solution to the optimal design \citep{harris1990many}. 

These three basic problem classes have been integrated in different combinations to address more complex problems faced in real-world logistics systems. Inventory operations at a facility are ultimately determined by the demand size and the service area of this facility, which is the outcome of location decisions. This connection is modeled with CA integrating both location and routing decisions \citep{rosenfield1992application}. An apparent tradeoff is that a higher investment of facilities usually reduces long-haul distances for delivery vehicles and thus decreases the total routing cost. In problems integrating routing and inventory decisions, CA relates service area sizes and frequencies of delivery trucks to inventory sizes and holding costs \citep{daganzo1988comparison}. The basic tradeoff is that a higher delivery frequency and a smaller service area often reduce inventory costs while increasing transportation costs. In problems where location, inventory and routing costs are all considered, these three cost components shall be integrated together for solving the optimal design for the entire system. Integrating these three cost components in the CA framework is simply a summation of their respective analytical cost functions, which oftentimes yields closed-form analytical solutions as well \citep{campbell1990designing,campbell93a}. The classic book by \citet{daganzo1996logistics} provides a comprehensive review of CA models in the context of one-to-one, one-to-many and many-to-many distribution systems. A later survey by \citet{langevin96} covers the history, basic concepts, and developments of the CA models for freight distribution problems up to the mid-1990s. 

Since 1996, CA methods have continued to undergo significant adaptations and advancements in the contexts of a variety of emerging problems. In this paper, we provide a review of the methodological advancements and applications of CA models in the period of 1996-2016, which were not covered in the previous surveys and yet had a major impact on the current state of the art. Most notably, a number of efforts have been made to address issues related to service reliability, competition, time-dependency and emerging services and technologies. We will discuss the areas that had been covered, and propose research directions that need further work. In paper organization, we classify the reviewed studies into three broad categories: (i) facility location studies, (ii) distribution and transit  studies, and (iii) integrated supply chain and logistics studies, as a large portion of recent studies in CA can be covered with one of these three classes. We note that this classification is not necessarily disjoint nor is it exhaustive. Similar to many other systems, the complexity of modern logistics and distribution systems is high. Many recent studies integrate various components of these systems. For example, some facility location studies consider routing and distribution costs as a part of the objective whereas some distribution problems consider different service regions as in facility location problems. Our classification is based on decision variables of the problem. If the decision variables represent the locations (areas) of the facilities to be located or to be strengthened, then the study is classified under facility location studies. If the decision variables represent the vehicle's routing, then the study is classified under transit and distribution studies. Studies that combine location decisions with routing or inventory decisions are reviewed under the integrated supply chain and logistics studies.

The remainder of this paper is organized as follows: The location models are reviewed in Section \ref{Sec_Location_Model}, and the distribution and transit models are reviewed in Section \ref{Sec_Routing_Model}. The integrated models are discussed in \ref{Sec_Hybrid_Model}. Finally, Section \ref{Sec_Conclusions} concludes this paper by discussing current trends and potential gaps in CA.

\section{Location Studies}
\label{Sec_Location_Model}
Facility location problems determine the optimal configuration of a facility network to meet a given objective (e.g., maximize service level, minimize cost). Such problems involve various types of cost components such as fixed facility opening and variable operating costs, inventory costs, transportation costs and so forth. In the discrete setting, the underlying network consists of demand (customer) and potential facility location points. The facility location problem and its extensions have been carefully investigated in reviews such as \cite{aikens1985facility}, \cite{owen1998strategic}, \cite{klose2005facility}, \cite{Shen07} and \cite{melo2009facility}. 

As the numbers of customers and potential facility locations increase, the number of variables in the underlying mathematical programming model increases, resulting in hard-to-solve large-scale problems. As mentioned before, CA enables us to handle these large-scale problems. In addition, decision makers must account for uncertainty in demand while locating the facilities in the network, see \cite{snyder2006facility} for a review on facility location under uncertainty. Aggregating customers and representing the total demand of customers in the aggregated region as a demand density per area, as is done in CA, can help reduce the uncertainty for the decision makers. For example, \cite{ouyang10cui} show that the CA method is a promising tool for finding near-optimum solutions, within {4\%} - {7\%} of optimality, especially when the demand distribution is highly variable across space.

Applications of CA in facility location problems have increased in two decades, and there are four streams of studies. The first stream focuses on the classical setting of facility location problems which is finding the location of facilities to be opened. The term classical is used here to describe the facility location problems that do not consider disruption of facilities in their models. In other words, the classical problems assume the settings are deterministic and known and the facilities never fail to serve their customers due to disruptions. As counterparts of these deterministic facility location problems, the second stream of studies considers stochasticity and uncertainty from internal operations and external environments by taking into account the disruption of facilities in modeling.  The third stream includes studies that adopted CA methods to address problems of siting facilities for a startup company against incumbent facilities in a competitive market. Finally, the fourth stream consists of studies that focus on the problem of region discretization which generally deals with solution implementation; i.e., transforming CA solutions to a discrete setting. Such studies allow for greater implementation of CA solution. 

After the review of basic concepts and formulations, this section reviews recent studies using CA models to obtain the optimal location of various facilities and its extensions. 

\subsection{Basic Concepts}

\label{Location_Basic_Models} 
\begin{figure}[htbp]
\centering{} \subfigure{ \label{UFL_1D_new} \includegraphics[width=0.45\textwidth]{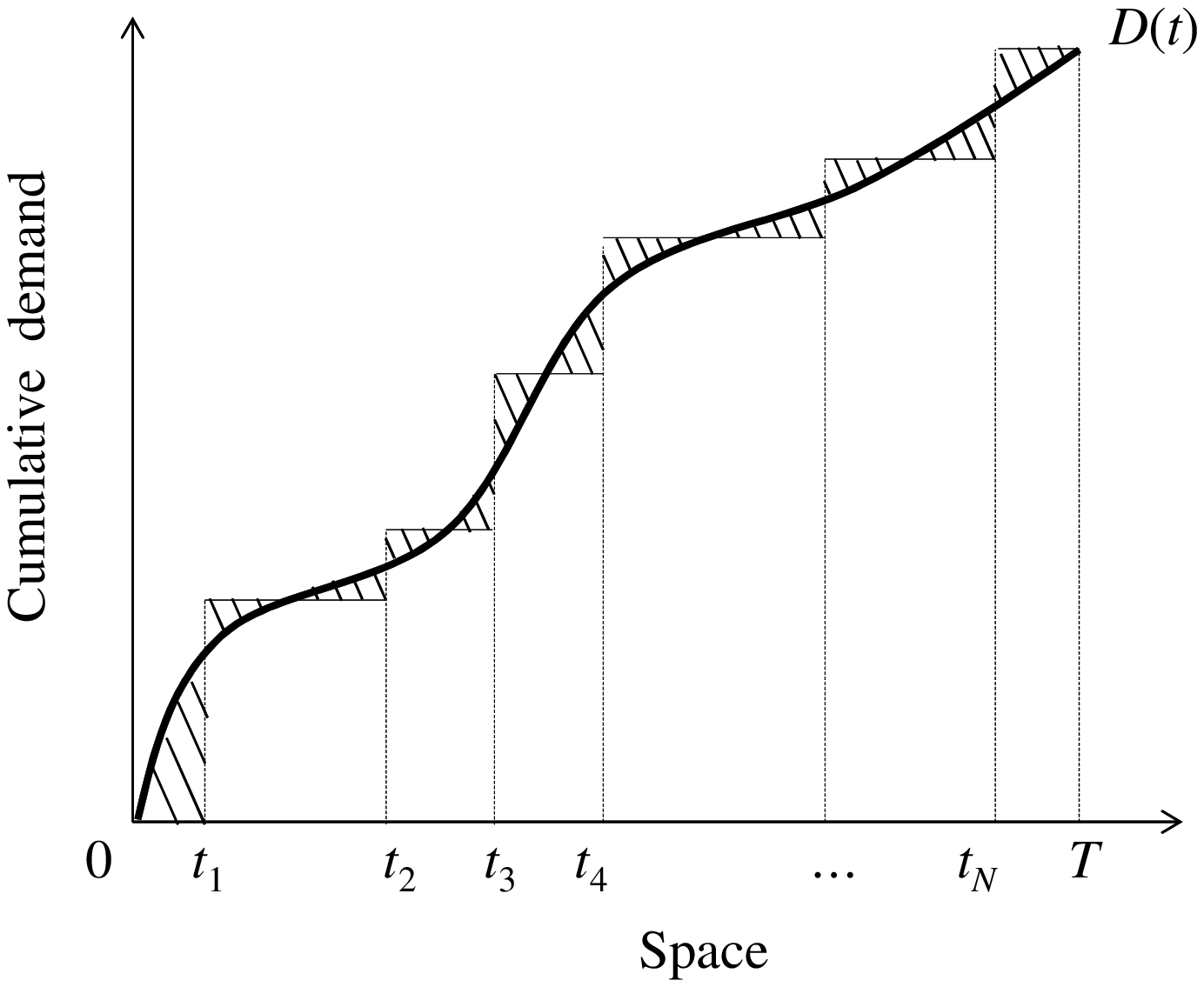}
} \subfigure{ \includegraphics[width=0.4\textwidth]{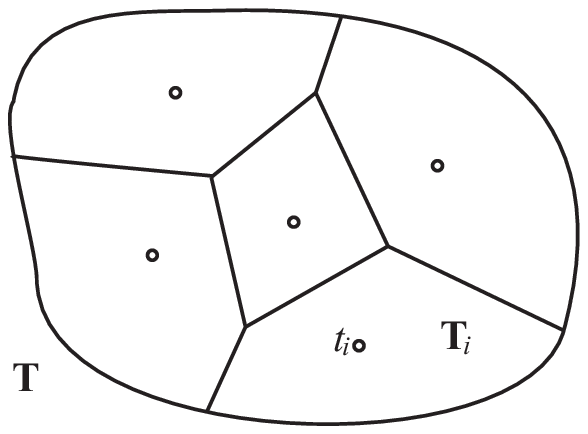} \label{UFL_2D_new}
} \caption{Uncapacitated facility location problems: \subref{UFL_1D_new} one-dimensional;
\subref{UFL_2D_new} two-dimensional.}
\end{figure}

%The above CA framework has been extended to solve the one-dimensionalfacility location problem For example, with slight modification, Figure \ref{Lot_Size} can be reinterpreted as a specific facility location problem in one-dimensional space $\mathbf{T}$.
CA is initially applied to one-dimensional Uncapacitated Facility Location (UFL) problems \citep{newell73,daganzo86}, as illustrated in Figure \ref{UFL_1D_new}. Customer demand density at each location $t\in\mathbf{T}:=[t_{0},T]$ is denoted by $d(t)$, and then its accumulation from locations $t_{0}$ to $t$ is $D(t)$ as shown in Figure \ref{UFL_1D_new}. Facilities are built at locations $t=t_{1},...,t_{N}$, each with opening cost $f(t)$. Once facilities are built, each customer will be served by its closest facility. It is assumed that the transportation cost to serve each unit demand at location $t$ is the product of its travel distance and a constant scalar $c$, i.e., $c|t-t_{i}|$ where $t_{i}$ is the closest point to $t$ among $\mathbf{T}$. Thus, the total travel cost can be represented as the shaded area in Figure \ref{UFL_1D_new}. For the convenience of the notation, define $t'_{i}:=(t_{i}+t_{i+1})/2,\forall i=1,\cdots,N-1$, $t'_{0}=t_{0}$ and $t'_{N}=T$. Then the objective of this problem is to determine the optimal facility locations that minimize the total cost integrating both facility opening expenses and transportation costs, as formulated below: 
\begin{equation}
\min_{t_{1}\le t_{2}\le\cdots\le t_{N}\in\mathbf{T}}\sum_{i=1}^{N}\left(f(t_{i})+c\int_{t'_{i-1}\le t\le t'_{i}}\left|t-t_{i}\right|d(t)dt\right)\label{1D_UFL_formulation}
\end{equation}

%Again, we can represent problem \eqref{1D_UFL_formulation} in the discrete form. Similar to model \eqref{LS_Model}, we discretize space $\mathbf{T}$ into $N$ locations $\{t'_1,t'_2,\cdots,t'_N\}$, and we redefine $c_{ij}$ as $c|t'_i-t'_j|, \forall i, j =1,2,\cdots,N$. The remaining notation keeps the same. Then the discrete version of \eqref{1D_UFL_formulation} can be formulated into \eqref{LS_Model} as well. On the opposite side, we can also formulate its CA version in a similar way of formulating \eqref{Continuous_Cost}.We
See from Figure \ref{UFL_1D_new} that the influence area of facility $i$ (i.e., the area that facility $i$ serves) is $\mathbf{T}_{i}:=[t'_{i-1},t'_{i})$, and its size is $A_{s}(t)=|t'_{i+1}-t'_{i}|,\forall i=1,\cdots,N$. Then, the total cost \eqref{1D_UFL_formulation} can be rewritten as: 
\begin{equation}
\min_{t_{1}\le t_{2}\le\cdots\le t_{N}\in\mathbf{T}}\sum_{i=1}^{N}\int_{\mathbf{T}_{i}}\left[\frac{f(t_{i})}{A_{s}(t)}+\frac{cA_{s}(t)}{4}d(t_{i}^{*})\right]dt,\label{C_i_new}
\end{equation}
where $t_{i}^{*}\in\mathbf{T}_{i}$ satisfies that $\frac{1}{4}(t'_{i}-t'_{i-1})^{2}d(t_{i}^{*})$ equals the shaded area in this interval $\mathbf{T}_{i}$ in Figure \ref{UFL_1D_new}.

We assume that functions $d(t)$ and $f(t)$ vary slowly across $\mathbf{T}$. Then without loss of accuracy, demand $d(t_{i}^{*})$ can be well approximated by $d(t)$ in \eqref{C_i_new}, and step function $A_{s}(t)$ can be replaced by a continuous service area function $A(t)$. This yields a CA model of equation \eqref{C_i_new} as follows: 
\begin{equation}
\min_{\{A(t)\}_{t\in\mathbf{T}}}\int_{\mathbf{T}}\left[\frac{f(t)}{A(t)}+\frac{cA(t)}{4}d(t)\right]dt.\label{Continuous_Cost_new}
\end{equation}
Model \eqref{Continuous_Cost_new} can be easily solved by optimizing each integrand component independently, which yields the optimal service area function as follows: 
\begin{equation}
A(t)=\left[4f(t)/(cd(t))\right]^{1/2},\forall t\in\mathbf{T}.\label{Opt_A_t_new}
\end{equation}

The modeling concepts underlying this one-dimensional problem have been generalized to UFL problems in two-dimensional spaces \citep{daganzo86} or even three dimensional spatio-temporal cubes \citep{daganzo05}. Figure \ref{UFL_2D_new} presents a two-dimensional example. In space $\mathbf{T}\in\mathbb{R}^{2}$, a set of facilities, each denoted by $t_{i}$, serves distributed customers that have density $d(t),\forall t\in\mathbf{T}$. Since each customer is served by its closest facility, $\mathbf{T}$ is tessellated into regions such that facility $t_{i}$ serves customers in region $\mathbf{T}_{i}$ (which is also called Voronoi Tessellation). All other parameters remain the same. The objective again is to find the optimal facility location to minimize the total system cost. By extending the CA scheme to two-dimensional space, we can formulate the continuous counterpart of this problem as follows: 
\begin{equation}
%\min_{\{A(t)\}_{t\in\mathbf{T}}}\int_{\mathbf{T}}\left[\frac{f(t)}{A(t)}+\frac{2c\sqrt{A(t)}}{3\sqrt{\pi}}d(t)\right]dt.
\min_{\{A(t)\}_{t\in\mathbf{T}}}\int_{\mathbf{T}}\left[\frac{f(t)}{A(t)}+c_{m}c\sqrt{A(t)}d(t)\right]dt,\label{Continuous_Cost_2D}
\end{equation}
where scalar $c_{m}$ represents the ratio of the average transportation distance in the service region, which is determined by the distance metric and the assumption of geometry of the influence area \citep{newell73}. With a Euclidean distance metric, $c_{m}$ is 0.376, 0.377, 0.382 and 0.403 respectively for a circular, hexagonal, square and triangular influence area. While circles yield the minimum average distance, they can not fill up the space. It has been proven that the hexagonal tessellation is the optimal feasible layout, when the space is homogeneous and large \citep{ouyang10cui}. Overall, these $c_{m}$ values are close (particular for circular, hexagonal and square areas). This suggests that the effect of influence area geometry is rather minor as long as the area does not become too elongated. If the Manhattan metric is used, then the optimal influence area shall be a square and the corresponding $c_{m}$ value is $0.454$. Again, the minimizer of \eqref{Continuous_Cost_2D} can be solved for each integrand as follows: 
\begin{equation}
A(t)=\left[\frac{f(t)}{2cc_{m}d(t)}\right]^{2/3},\forall t\in\mathbf{T}.\label{Opt_A_t_2d}
\end{equation}

\subsection{Classical Facility Location}

\cite{langevin96} review the early studies that analyze the location of facilities using continuous models \citep{beckmann1968location, larson1972insensitivities, love1972computational, newell73, newell1980traffic, erlenkotter1989general, klincewicz1990fleet, geoffrion1995distribution}. In this section, we review recent papers that study the classical (deterministic) facility location problem using CA. Deterministic models assume that customer demand is deterministic and known and facilities are always properly functioning. Please see Table \ref{tab:classicalFacility} for the summary of the papers. As with all tables in the paper, we rename some variable and cost terms to obtain a common framework to compare across papers.

\setlength{\extrarowheight}{1em}
\begin{table}[!ht]
\caption{Classical facility location papers summary.}
\vspace{0.25cm}
  \centering
  \footnotesize
    \begin{tabular}{c| L{4cm}| P{0.25cm} P{0.25cm} P{0.25cm} P{0.25cm} P{0.25cm}|P{7cm}}
          & & \begin{sideways} \hspace{-0.5cm} Facility opening  \end{sideways} & \begin{sideways} \hspace{-0.5cm} Facility operating \end{sideways} & \begin{sideways} \hspace{-0.5cm} Transportation \end{sideways} & \begin{sideways} \hspace{-0.5cm} Inventory \end{sideways} & \begin{sideways} \hspace{-0.5cm} Other\end{sideways} &  \\
   	&\multicolumn{1}{c|}{\textbf{Study}} & \multicolumn{5}{c|}{\textbf{Cost Structure}} & \textbf{Main Contribution(s)} \\
    \toprule
	\multirow{5}{*}[-2em]{\begin{sideways}Location\end{sideways}}    
    & Dasci and Verter (2001) 	& X     & X     & X     &       &       & \parbox{7cm}{\linespread{0}\selectfont Varies distance metric.}  \\
    & Wiles and Brunt (2001) 	&       & X     & X     &       &       &  \parbox{7cm}{\linespread{0}\selectfont Introduces circular continuum region and direct route between the continuum location and the transshipment depot.} \\    
    & Dasci and Verter (2005) 	& X     &       & X     &       &       & \parbox{7cm}{\linespread{0}\selectfont Analyzes product focus and market focus strategies.} \\
    & Carlsson and Jia (2013a) 	& X     &       & X     &       & X     & \parbox{7cm}{\linespread{0}\selectfont Introduces backbone network costs and fast constant-factor approximation algorithm for placing facilities in any convex polygonal region.}  \\
    & Carlsson and Jia (2013b) 	& X     &       & X     &       & X     & \parbox{7cm}{\linespread{0}\selectfont Considers asymptotic behavior of the optimal solution and does input parameters analysis.}  \\
	\midrule   
	\multirow{4}{*}[0em]{\begin{sideways}\shortstack{Location \& \\Allocation}\end{sideways}} 
    & Murat et al. (2010) 		& X     &       & X     &       &       & \parbox{7cm}{\linespread{0}\selectfont Considers that allocation decisions are given before location decisions.} \\
    & Hong et al. (2012) 		& X     &       & X     &       &       & \parbox{7cm}{\linespread{0}\selectfont Introduces Robust Integer Facility Location (RIFL) and Robust Continuous Facility Location (RCFL) models.}  \\
    & Bouchery and Fransoo (2015) &       &       & X     &       &       & \parbox{7cm}{\linespread{0}\selectfont Analyzes intermodal network design decisions from a cost, carbon emissions and modal shift perspective.} \\
    \midrule
    \multirow{4}{*}[0em]{\begin{sideways} \shortstack{Other Operational \\ Constraints} \end{sideways}} 
    & Cachon (2014)				& X     & X     & X     & X     &       & \parbox{7cm}{\linespread{0}\selectfont Combines an inventory model with the TSP and the k-median problem.} \\
    & Wang et al. (2014) 		& X     & X     & X     &       &       & \parbox{7cm}{\linespread{0}\selectfont Adds a time dimension.} \\
   & Xie and Ouyang (2015)				& X     & X     & X     & X     &       & \parbox{7cm}{\linespread{0}\selectfont Proves optimal layout of a transshipment system that combines the TSP and the k-median problem, with inventory considerations.} \\
    & Ouyang et al. (2015)		& X     &       & X     &       & X     & \parbox{7cm}{\linespread{0}\selectfont Integrates both facility location decisions and traffic equilibrium in a continuous space.} \\
	\bottomrule
    \end{tabular}%
    \label{tab:classicalFacility}
\end{table}%

There are three main divisions of studies. The first division includes papers that solely focus on facility location problems. In most studies using CA, cost metrics are assumed to be constant over the entire region. Thus, the corresponding CA model returns the service regions for each facility and a Weber-type problem \citep{weber1929theory} is solved to finalize their location decisions. \cite{dasci01} determine the locations of the facilities that minimize the total fixed facility opening costs, facility operating costs and transportation costs in production distribution systems. To find the average distance from the facility to a customer in the service region and therefore the transportation cost, unlike most studies that use the exact locations of the facilities,  they use a distance metric that varies over the service region. In a more generalized case of the continuous facility location problem, \cite{Dasci05_Eval} develop a methodology based on a facility design model to select between product focus \footnote{A firm's ability to increase productivity and lower cost by limiting the number and variety of operations at its production lines.} and market focus strategies. They use continuous functions to represent the spatial distribution of demand and cost parameters. \cite{wiles2001optimal} formulate a problem of finding an optimal location of transshipment depots within a harvesting region. They assume that the agriculture commodities are distributed in a circular continuum region and that a transshipment depot would be located for the collection and export of the product. Instead of taking the radial and circumferential route \citep{lam1967flow, blumenfeld1970routing} for the circular modeling region, the authors consider a direct route between the continuum harvesting location and the transshipment depot. 

Considering a constant cost metric over the service region, \cite{carlsson2014continuous} locate facilities to minimize the total facility opening, transportation and backbone network costs. The key difference of this paper with the literature is the idea of considering the backbone network costs, costs from connecting the facilities to each other. They provide a fast constant-factor approximation algorithm for placing facilities in any convex polygonal region. Implementing the same idea, \cite{carlsson2013euclidean} consider the problem of designing an optimal spoke-hub distribution network\footnote{The spoke-hub distribution network is a system of connections arranged like a chariot wheel, in which all traffic moves along spokes connected to the hub at the center.} . They give an approximation algorithm to select the optimal locations of a set of hubs in a convex planar Euclidean region where the “spokes” are continuously and uniformly distributed in the region. Their work is different from location-routing problems where the network of customers is served by vehicle tours since they think of the backbone network of facilities instead of customers. 

The second set of papers not only studies the location problem, but also considers allocation decisions which determine the assignment of customers to facilities once facility locations are known. \cite{Murat10} find the locations and allocations of facilities over a region with dense demand. The novelty of this paper is that they prioritize the determination of the service regions over the location decisions. The authors determine the boundaries of each service region using a gradient algorithm. For emergency response facility and transportation problems, \cite{hong2012development} compare the performance of Robust Integer Facility Location (RIFL) where the demand of a commodity distribution point (CDP) is covered by a main distribution warehouse or a backup CDP and the Robust Continuous Facility Location (RCFL) where that of a CDP is covered by multiple distribution warehouses models. They find that, under normal conditions, RIFL performs better while, under the emergency conditions, the RCFL outperforms others and the overall logistics cost and the robustness level of the RCFL are better. Another recent location-allocation study is \cite{bouchery2015cost}, who develop a model to optimize the terminal location and the allocation between direct truck transportation and inter-modal transportation. Approximating the demand as continuous, they analyze the dynamics of intermodal freight transportation with respect to cost, modal shift and carbon emissions.

The third set of papers extends the facility location problem by adding new operational constraints. Adding a time dimension, \cite{wang2015infrastructure} formulate a continuous model for the Dynamic Facility Location Problem (DFLP). Their proposed model finds the best location and timing of facility deployment that minimizes the logistics cost during a planning horizon. Unlike \cite{campbell1990locating} who considers mobile terminal locations, they force the locations of open facilities to remain unchanged in the remainder of the planning horizon (location consistency constraint). In their study, CA is used to find the optimal facility density in the spatiotemporal continuum. To discretize the CA output into a set of discrete facility locations, they extend the disk model (for one static time period) of \cite{ouyang06} to a tube model (for multiple time periods). Analyzing the accuracy and convergence of the proposed generic and flexible method, they show that it solves the DFLP to a reasonable accuracy. 

The number of facilities to locate may also be an endogenous decision variable related to other factors such as the routing of inbound supply. \cite{cachon2014retail} uses the CA approach to locate retail stores and evaluate their environmental impacts under different governmental policies. The retail store density problem combines an inventory model with the TSP and the continuous k-median problem. \cite{xie2015optimal}  further prove optimal spatial layout of transshipment facilities and the corresponding service regions in an infinite homogeneous plane to minimize the total cost for facility set-up, outbound delivery and inbound replenishment transportation. On a Euclidean plane, a tight upper bound (within 0.3\%) can be achieved by a type of elongated cyclic hexagons. A similar elongated non-cyclic hexagon shape, with proper orientation, is optimal for service regions on a rectilinear metric plane.  \cite{ouyang2015facility} study a median-type facility location problem under elastic customer demand and traffic equilibrium in a continuous space using both CA and mixed-integer programming. The median-type facility location problem seeks optimal facility locations to minimize the transportation costs for serving spatially distributed customers, subject to a budget constraint for facility investments. Unlike the traditional models, they model the destination and routing decisions of customers endogenously with the decisions on facility locations. The locations of facilities are not given in advance in their model. 

In sum, we observe a shift from problems focusing solely on the locations of the facilities to more integrated models that not only provide solutions to the facility location problem but also consider operational constraints such as allocation decisions, number of facilities, time windows and so on. The application of CA, due to its modeling convenience, is more vivid, when the mode is more complex.

\subsection{Disruption in Facility Reliability}

%{\color{blue} Will be taken from Section 3.2 of Smilowitz Draft. and Section 2.3 of Ouyang Draft.}

In this section, we review facility location problems where facilities are subject to disruptions. When a facility fails to provide service for some period of time, the customers in the corresponding region are either served by the closest facilities or not served at all. Reliability models have been developed to include the effect of random disruptions when determining facility locations. See Table \ref{tab:DisruptionFacility} for the summary of reliability papers.

\setlength{\extrarowheight}{1em}
\begin{table}[ht]
\caption{Disruption in reliability papers summary.}
\vspace{0.25cm}
  \centering
  \footnotesize
    \begin{tabular}{c| L{4cm}| P{0.25cm} P{0.25cm} P{0.25cm} P{0.25cm} P{0.25cm}|P{7cm}}
          & & \begin{sideways} \hspace{-0.5cm} Facility opening  \end{sideways} & \begin{sideways} \hspace{-0.5cm} Facility operating \end{sideways} & \begin{sideways} \hspace{-0.5cm} Transportation \end{sideways} & \begin{sideways} \hspace{-0.5cm} Inventory \end{sideways} & \begin{sideways} \hspace{-0.5cm} Other\end{sideways} &  \\
    & \multicolumn{1}{c|}{\textbf{Study}} & \multicolumn{5}{c|}{\textbf{Cost Structure}} & \multicolumn{1}{c}{\textbf{Main Contribution(s)}} \\
    \toprule
    \multirow{5}{*}[-1em]{\begin{sideways}Locations\end{sideways}}
    & Dasci and Laporte (2005) 		& X     & X      & X     &       &       & \parbox{7cm}{\linespread{0}\selectfont Considers stochastic demand and capacity acquisition as a decision variable.} \\
    & Cui et al (2010) 		& X     &       & X     &       &       & \parbox{7cm}{\linespread{0}\selectfont Considers region dependent disruption probabilities.} \\
    & Li and Ouyang (2010) 	&       &       & X     &       & X     & \parbox{7cm}{\linespread{0}\selectfont Considers probabilistic facility disruption risks.} \\
    & Li and Ouyang (2012) 	&       & X     &       &       &       & \parbox{7cm}{\linespread{0}\selectfont Considers location dependent failure probabilities.} \\
    & Wang and Ouyang (2012)& X     &       & X     &       &       & \parbox{7cm}{\linespread{0}\selectfont Uses game-theoretical models under spatial competition and considers probabilistic facility disruption risks.} \\
    & Lim et al. (2013)		& X     &       & X     &       &       & \parbox{7cm}{\linespread{0}\selectfont Assumes disruption probabilities known imperfectly and computes the cost of underestimation and overestimation.} \\
	\midrule   
	\multirow{2}{*}[-1em]{\begin{sideways}Other\end{sideways}} 
    & Bai et al. (2015) 	& X     &       & X     &       & X     & \parbox{7cm}{\linespread{0}\selectfont Models facility disruptions in biofuel supply chain design.} \\
    & Wang et al. (2015)	& X     & X     & X     &       &       & \parbox{7cm}{\linespread{0}\selectfont Considers spatially correlated, site-dependent probabilistic distribution for the resource supply.} \\
	\bottomrule
    \end{tabular}%
    \label{tab:DisruptionFacility}
\end{table}%

Earlier studies focused on the design of systems with sufficient redundancy under stochastic demand fluctuations. \cite{dasci2005analytical} consider the location and capacity acquisition problem with stochastic demand. They determine the optimal market area of the facilities instead of exact locations. They also consider two different strategies for insufficient capacity: outsourcing the extra demand or considering the extra demand as lost sales. The optimal market areas are found using CA which minimizes total cost. This cost includes the cost of the facility, capacity acquisition and operation, transportation, and shortages. They find that the solution to the outsourcing strategy requires more facilities. Their model also shows that the out-sourcing solution depends on demand, fixed costs, and transportation cost, while the lost sales solution additionally depends on the distribution of demand and the shortage cost.

\cite{li2011reliable} locate facilities in an environment with correlated disruptions. Assuming that the facility locations are fixed and facilities are uncapacitated, the authors develop a CA formulation to solve the problem over an infinite homogeneous plane with independent disruptions. They use the results of this simplistic model as a building block for correlated disruptions. In the non-homogeneous and correlated case, the CA method returns the service regions for each facility and the authors suggest the use of discretization method developed in the study of \cite{ouyang06}. \cite{ouyang10cui} locate facilities under region dependent random disruptions. Assuming that the facilities are uncapacitated and the distance metric is Euclidean, the authors develop both discrete and continuous models for the same problem. Similar to \cite{li2011reliable}, the homogeneous case leads to a simpler continuous model. The solution of this model leads to the hexagonal tessellation of the entire region due to Euclidean metric. For the non-homogeneous case, they calculate optimal service regions using the discretization algorithm of \cite{ouyang06}. As expected, CA approach performs better in larger instances when compared to the corresponding discrete model. \cite{li2012reliable} maximize the surveillance of traffic flow by locating sensors (facilities) with a limited budget. Similar to \cite{ouyang10cui}, the sensors are subject to location dependent random (but uncorrelated) disruptions. They show that, under certain conditions, CA provides a lower bound for MILP. Considering discrete and continuous reliable facility location models, \cite{bai2015effects} design a model for reliable bio-ethanol supply chains with potential operational disruptions. Their model includes the possibility of simultaneous disruptions of multiple facilities serving the same customer. The impacts of both site-independent (similar to \citet{ouyang10cui}) and dependent disruptions are analyzed in an empirical case study. In the prior papers, the disruption probabilities are known. However, \cite{lim2013} determine the locations of reliable and unreliable facilities in an environment with unknown disruption probabilities. Under such a setting, they formulate a CA model and show that underestimating both disruption probability and correlation degree result in a greater increase in the expected total cost compared to the overestimating.

%%%MERGE THE FOLLOWING INTO THE NEXT SECTION  - YO???
%Taking competition in to account, \cite{wang2013continuum} model a facility location problem in a competitive environment with disruptions. In addition to the Stackelberg setting (Leader-Follower) considered in \cite{Dasci2005}, they consider Nash setting for the competing companies. Assuming that the customers always benefit from seeking service (i.e. they prefer service no matter how far the closest store is rather than having no service at all), the authors provide a game theoretical models in the CA scheme under spatial competition. Incorporating uncertainties in supply/demand and the risk of facility disruptions, \cite{wang2015infrastructure} develop a game-theoretic modeling framework using CA to study the impacts of competition on the optimal infrastructure deployment. Similar to \cite{carlsson2013euclidean} and \cite{carlsson2014continuous}, they choose bio-fuel industry to conduct their case study and focus on the resource supply competition in this industry. Unlike \cite{wang2004locating}, which study deterministic market share competition between two firms, they focus on the post-entry facility deployment plan for an emerging industry under uncertainties. Specifically, they assume a spatially correlated, site-dependent probabilistic distribution for the resource supply, while each built facility may be disrupted independently with a site-dependent probability. 

In sum, CA has been widely used to facilitate developing facility location models with more realistic disruption setting implemented in novel application areas.  

\subsection{Competitive Facility Location Models}

%{\color{blue} Will be taken from Section 2.4 of Ouyang Draft and some papers of Section 3.1 of Smilowitz Draft.}

In this section, we review recent papers that study the competitive facility location problem using CA. Please see Table \ref{tab:competativeFacility} for the summary of the papers. The location models introduced in Section \ref{Sec_Location_Model} assume a centralized system structure; i.e., a centralized agency makes all location decisions and its goal is to maximize a single system benefit measurement (or minimize a single system cost measurement). While these models can address location planning for a public agency or a monopolistic private corporation, they are not suitable for an open market environment where multiple players are competing against each other and each makes separate facility planning decisions so as to maximize its own profit. Recently, several multi-store competitive location models have been proposed in literature to address the facility location design problem in such a competitive environment.

\setlength{\extrarowheight}{1em}
\begin{table}[ht]
\caption{Competitive facility location papers summary.}
\vspace{0.25cm}
  \centering
  \footnotesize
    \begin{tabular}{ L{4cm}| P{0.25cm} P{0.25cm} P{0.25cm} P{0.25cm} P{0.25cm}|P{7cm}}
           & \begin{sideways} \hspace{-0.5cm} Facility opening  \end{sideways} & \begin{sideways} \hspace{-0.5cm} Facility operating \end{sideways} & \begin{sideways} \hspace{-0.5cm} Transportation \end{sideways} & \begin{sideways} \hspace{-0.5cm} Inventory \end{sideways} & \begin{sideways} \hspace{-0.5cm} Other\end{sideways} &  \\
   	\multicolumn{1}{c|}{\textbf{Study}} & \multicolumn{5}{c|}{\textbf{Cost Structure}} & \textbf{Main Contribution(s)} \\
    \toprule
     Dasci and Laporte (2005) 	& X     &       &       &       & X     & \parbox{7cm}{\linespread{0}\selectfont Considers a competitive environment.} \\
     Wang and Ouyang (2013) 	& X     &       & X      &       &      & \parbox{7cm}{\linespread{0}\selectfont Introduces reliable competitive facility location problem.} \\
     Wang et al. (2013) 	&      & X      &       & X      & X     & \parbox{7cm}{\linespread{0}\selectfont Considers a non-cooperative game in a multi-echelon bio-fuel supply chain network.} \\
	\bottomrule
    \end{tabular}%
    \label{tab:competativeFacility}
\end{table}%

As a pioneer in using CA in a competitive environment, \cite{Dasci2005} study the multi-store competitive location problem. Treating the consumers as entities continuously spread over the market, they consider competition based on location density. Assuming that the customers go to the nearest store, the paper provides insights for different strategies of leader and follower companies. Their model decides whether or not a company should enter the market as well as the optimal number and locations of stores.

Considering facility disruptions, \cite{wang2013continuum} find the optimal facility location design under spatial competitions  and general transportation cost functions. This problem involves discrete bi-level optimization and probabilistic facility disruption considerations,
both of which are notoriously difficult. Following the CA approach, the paper replaces discrete location variables by continuous and differentiable density functions to approximate the bi-level competition problem with closed-form formulas. Results of a case study for competitive bio-fuel supply chain design show insights on how competing bio-fuel companies should optimally plan their refinery location decisions.
% Even for the discrete bilevel problem itself, the discrete and combinatorial nature of the lower-level problem makes the problem very challenging to handle especially when the problem instance becomes large. The incorporation of a large number of probabilistic facility failure scenarios makes the problem even harder; even for a given design, to evaluate the system performance across all possible disruption scenario requires an exponential amount of time.
Extending competitive location models in bio-fuel industry to a multi-echelon supply chain, \cite{wang2013food} propose a CA model to optimize a supply chain network where resource suppliers independently choose among multiple competitive outlets. The work is motivated by how the diversion of agricultural crop land to energy feedstock production affects equilibria in both food and energy markets, and how the government could mitigate negative environmental impacts by the conservation land reserve program. The location decision of the supply chain is represented by a spatial density function, and the probability for a generic farmer to choose a particular land use outlet (food crop, energy crop, or reservation), based on price and distance, is derived in closed-form as a function of the local refinery density. Then, the optimal decision strategies of all stakeholders are incorporated into a tri-level non-cooperative Stackelberg game-theoretic model. %In the bio-energy context, the decisions of the raw material suppliers (i.e., farmers) include land allocation for food/energy productions or land reservation, and sales to nearby food markets or bio-refineries; the decisions of the supply chain include the number and locations of bio-refineries, and resource procurement prices; the government determines the optimal price for land reservation. The optimal spatial land use, the corresponding market equilibria, and the social welfare can also be derived.

% \hl {[**** REFERENCE wang2013optimal IS WRONG -- Should be a working paper: Wang, X., Lim, M.K. and Ouyang, Y. “Food, energy and environmental trilemma: Sustainable farmland use and biofuel industry development.” Energy Economics. Under revision. ]}

In sum, there are just a few studies that use the CA for modeling competitive environments in the literature. However, recent advances in game theory and interests in addressing problems in an open market environment may encourage more researchers to tackle problems existing in this area.

\subsection{Region Discretization}

%{\color{blue} Will be taken from Section 3.3 of Smilowitz Draft.}

To find the optimal service area for each facility in facility location problems using the CA approach, the discretization of the service region is required, which is not trivial. Discretization of the service region is to partition a continuous solution (e.g., continuous facility density) to discrete locations. The problem of optimally partitioning a region into smaller units (districts or zones), subject to constraints such as balance, contiguity and compactness is called the districting problem \citep{bunge1966theoretical}. The continuous districting problem is the problem in which the underlying space for facility sites and demand points are determined by continuous variables. Voronoi diagrams in combination with CA models can be useful in solving continuous location-districting problems; see \cite{aurenhammer1991voronoi} for a review on Voronoi Diagrams. Given a set of finite distinct points in a continuous planar space, associating all locations in that space with the closest member(s) of the point set with respect to the Euclidean distance, results in a tessellation of the space into a set of regions. This tessellation is called an ordinary Voronoi diagram. (For further explanation please refer to \cite{Okabe92}, p. 66).)

In this section, we review papers that allow for greater implementation of CA solutions by providing tools and methods for service region partitioning. See Table \ref{tab:RegionDisc} for the summary of region discretization papers.

\setlength{\extrarowheight}{1em}
\begin{table}[!ht]
\caption{Region discretization papers summary.}
\vspace{0.25cm}
  \centering
  \footnotesize
\begin{tabular}{L{4cm}| P{0.25cm} P{0.25cm} P{0.25cm} P{0.25cm} P{0.25cm}|P{7cm}}
          & \begin{sideways} \hspace{-0.5cm} Facility opening \end{sideways} & \begin{sideways} \hspace{-0.5cm} Facility operating \end{sideways} & \begin{sideways} \hspace{-0.5cm} Transportation \end{sideways} & \begin{sideways} \hspace{-0.5cm} Inventory \end{sideways} & \begin{sideways} \hspace{-0.5cm} Other\end{sideways} &  \\
    \multicolumn{1}{c|}{\textbf{Study}} & \multicolumn{5}{c|}{\textbf{Cost Structure}} & \multicolumn{1}{c}{\textbf{Main Contribution(s)}} \\
    \toprule
    Okabe and Suzuki (1997) &  \multicolumn{5}{c|}{\rule{2cm}{0.4pt}}     & \parbox{7cm}{\linespread{0}\selectfont Reviews continuous locational optimization problems with Voronoi diagram.} \\
    Du et al. (1999) 		&     \multicolumn{5}{c|}{\rule{2cm}{0.4pt}}    & \parbox{7cm}{\linespread{0}\selectfont Discusses Central Voronoi Tessellations (CVT).} \\
    Galvao et al. (2006) 	&       &       & X     &       &       & \parbox{7cm}{Develops a new partitioning scheme based on multiplicatively-weighted Voronoi diagrams.} \\
    Ouyang and Daganzo (2006) &    \multicolumn{5}{c|}{\rule{2cm}{0.4pt}}      & \parbox{7cm}{\linespread{0}\selectfont Improves CVT discretization. Generalizes to L1 norm with square influence regions.} \\
    Novaes et al. (2009) 	&       &       &       &       & X     & \parbox{7cm}{\linespread{0}\selectfont Uses the plane-sweep and the quad-tree techniques to construct Voronoi diagram.} \\
    Novaes et al. (2010) 	&       &       &       &       & X     & \parbox{7cm}{\linespread{0}\selectfont Gives an application of the Power Voronoi diagram.} \\
	\bottomrule    
    \end{tabular}
    \label{tab:RegionDisc}
\end{table}

In a review paper, \cite{okabe1997locational} identify eight classes of continuous locational optimization problems that can be solved through the Voronoi diagram method. \cite{Du99} also provide a general discussion on Central Voronoi Tessellations (CVT). CVT divides a region into sub-regions where the points in each sub-region are closer to this sub-region center compared to the centers of the other regions. This provides useful insights about one of the most frequently used tessellations in CA literature. One of the potential application areas that the authors mention is the mailbox placement problem, which belongs to the general class of facility location problems.

\cite{ouyang06} introduce a procedure that is an improvement on the well-known CVT method, avoiding undesirable shape formation. Assuming Euclidean distance metric, the authors define the optimal service regions to be circular with facilities at the center points. Using this property, they develop a method which uses different types of forces in the center and on the boundaries of the service regions to move and iterate through the regions. This study allows for greater implementation of a CA solution. Focusing on logistic districting problems, \cite{galvao2006multiplicatively} develop a new partitioning scheme based on multiplicatively-weighted Voronoi diagrams. The method uses predetermined weights and compares the closeness of points to the centers by considering weighted distances. The procedure is repeated until load factors in different regions are balanced sufficiently. The study compares the solutions generated by using Voronoi diagrams with a preliminary solution based on radial-ring. It is shown that the iterative Voronoi diagram generation procedure balances the load factors better than the ring-radial based procedure; however, the total travel distance is slightly lower in the latter procedure. The proposed approach can be used in environments with different road types and physical obstacles.

Combining a non-ordinary Voronoi diagram approach with an optimization algorithm, \cite{Novaes09} develop two continuous location-districting models applied to transportation and logistics problems. They use Voronoi diagram construction methods such as the plane-sweep and the quad-tree techniques in transportation and logistics problems which were previously applied in non-logistic problems such as computer graphics and robotics. They also extend the Voronoi diagram methodology to solve logistics districting problems with spatial barriers. Later, \cite{novaes2010continuous} investigate the application of the power Voronoi diagram in logistics districting problems. Associated with a continuous demand approach, they allow physical barriers such as rivers, reservoirs, hills, etc. into the vehicle displacement representation. The power Voronoi diagram is in the class of non-ordinary Voronoi diagrams useful in solving districting problems with barriers since the resulting Voronoi polygons are always convex. They compare the resulting district contours to the traditional wedge-shape formulation and show that not only they are smooth and closer to the configuration contours in practical situations, but also the resulting partition of the region leads to more balanced time/capacity utilization (load factors) across the districts. In conclusion, even though the reviewed papers provide useful tool for discretization of the service region, it is necessary to develop methods to be used for more complex CA models.

\section{Distribution and Transit Studies}
\label{Sec_Routing_Model}

This section reviews recent developments in using CA methods to solve vehicle routing related problems. The basic continuum models for routing problems are introduced in Section \ref{Routing_Basic_Model}. The developments since 1996 based on these models are reviewed in the following sections. Section \ref{sec:distribution_studies} discusses studies with a distribution focus where the aim is to determine the least cost routes to visit the demand nodes in a given network with a certain number of vehicles. For most of the realistic distribution systems, these problems are not tractable due to the difficulty of the underlying model or the size of the instance and exact methods are likely to fail in providing high quality solutions within an acceptable amount of time. Section \ref{sec:transit_studies} discusses studies with a transit focus. In transit studies, the aim is to design efficient networks by using various types of service options and routing the vehicles in the existing systems more efficiently. Overall, newly developed deterministic and stochastic CA models have been applied to a wide range or routing problems, including both supply chain management in private sectors and transit route design for public agencies. Various models have been built for planning and operation applications of different scopes. Discretization methods have also been proposed to convert a CA solution into implementable discrete route design.

\subsection{Basic Concepts}
\label{Routing_Basic_Model}
\begin{figure}[htbp]
\centering{}\subfigure{ \includegraphics[width=0.4\textwidth]{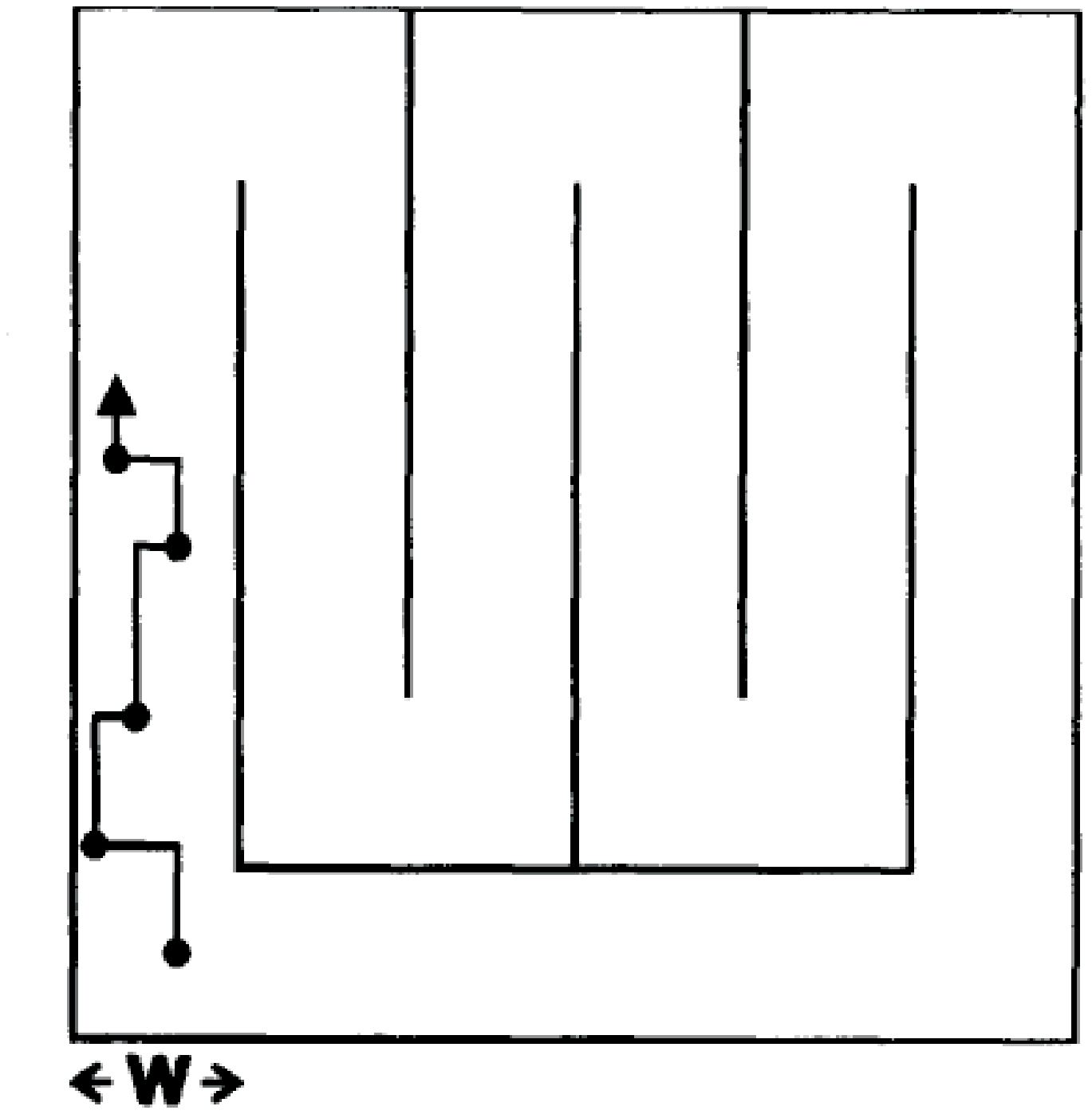}
\label{Strip_Strategy} } \subfigure{ \includegraphics[width=0.43\textwidth]{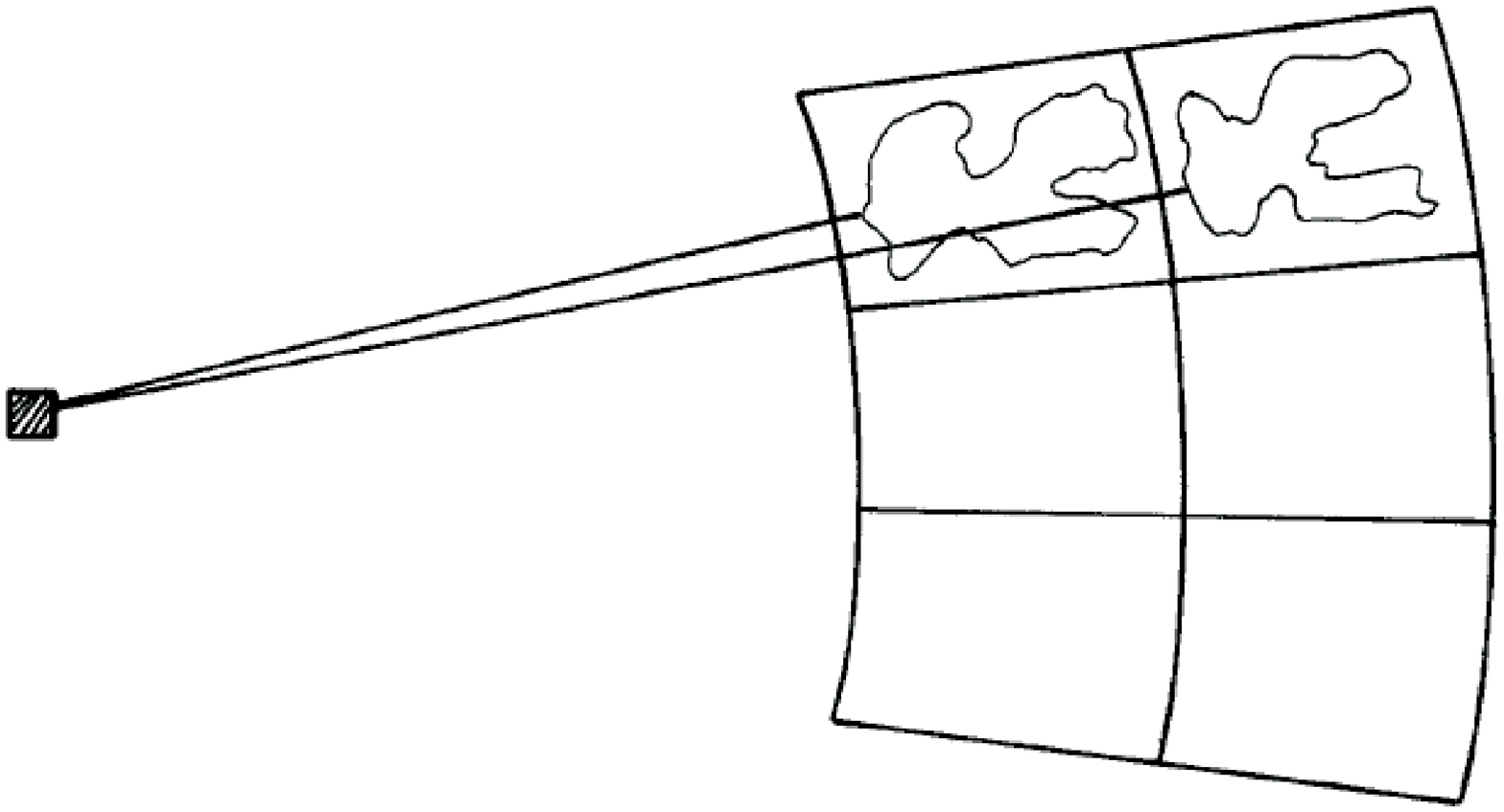}
\label{Ring_Sweep_Strategy}} \caption{(a) The strip strategy in \citet{daganzo84b} and (b) the ring sweep strategy in \citet{daganzo84a}.}
\end{figure}

Routing models are built on traditional routing problems, including Traveling Salesman Problem (TSP) and Chinese Postman Problem (CPP), and further consider a number of new components and aspects, such as special spatial partitions, customer service levels, multiple transportation modes, and service time windows. In TSP, $N$ customers are randomly distributed in a region $\mathbf{T}$ with a size of $T:=|\mathbf{T}|$, and a salesman needs to visit each customer within one tour and come back to where he starts. The objective is to minimize the salesman's total travel distance in this tour. Asymptotic analysis has shown that the optimal solution of the optimal total travel distance can be approximated with: 
\begin{equation}
k_{TSP}\sqrt{TN},\label{TSP_Dis}
\end{equation}
where $k_{TSP}$ is a constant, whose value is determined by the distance metric \citep{Beardwood1959,daganzo05}. Such analysis, however, did not indicate how to construct such an optimal tour, and the asymptotic value of $k_{TSP}$ is likely to underestimate the length of this tour for a finite-size problem with an irregular-shaped region. \citet{daganzo84b} proposed a strip strategy (see Figure \ref{Strip_Strategy}) to construct a near-optimal tour in a generic irregular region. The expected tour length using the strip strategy is again in the form of \eqref{TSP_Dis} where the $k_{TSP}$ value is only slightly higher than its asymptotic optimal value for different distance metrics and region shapes.

\begin{figure}[htbp]
\begin{center}
	\includegraphics[width=0.80 \textwidth]{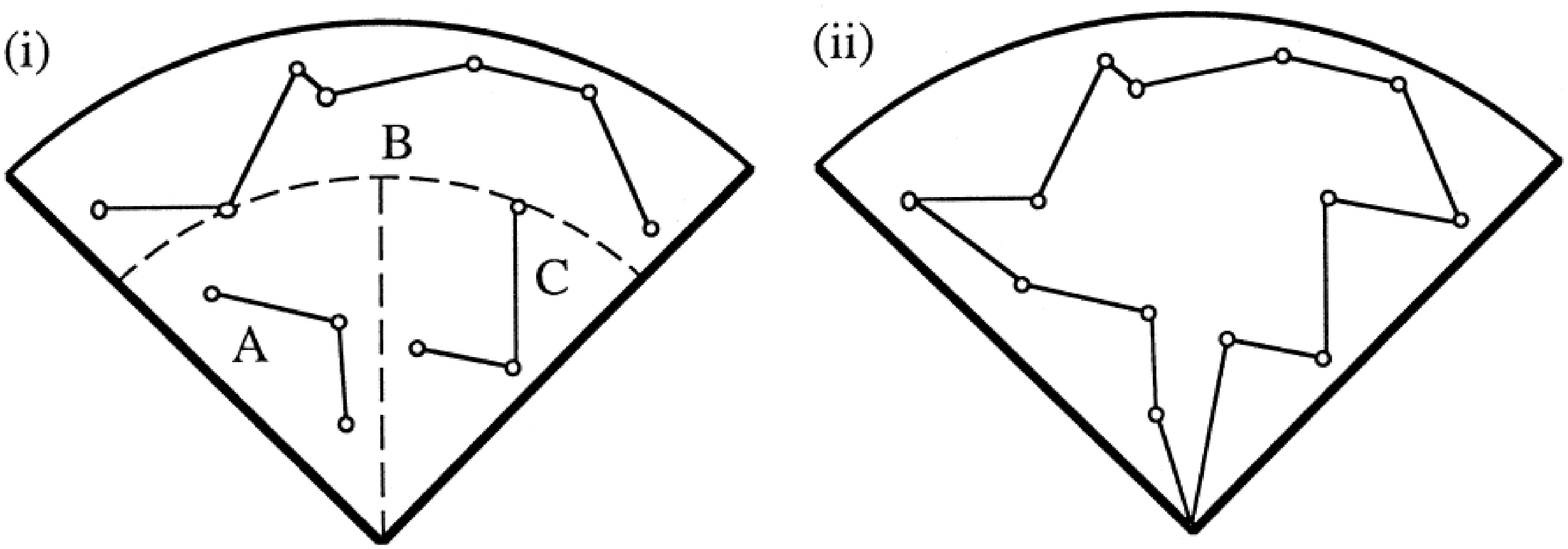} %width=0.80 \textwidth
	\caption{Tour construction strategy for a ring sector in \citet{del1999heuristic}.}\label{Castillo}
  \end{center}    
\end{figure}

\citet{del1999heuristic} proposed a new TSP tour construction strategy for uniformly distributed locations in a circular sector area, as illustrated in Figure \ref{Castillo}. This strategy splits a circular sector into three parts: two symmetric circular sub-sectors near the origin and a ring sub-sector around them.  The tour visits locations in these three sub-sectors sequentially without backtracking. It starts from the circular sub-sector on the left, goes through the above ring sub-sector, and comes back along the circular sub-sector on the right. Under this tour construction strategy, this paper derived analytical formulas of tour lengths for both Euclidean and ring-radial distance metrics. From experiments, this paper found that a tour length constructed with this strategy is comparable to that from equation \ref{TSP_Dis} with the strip strategy in \cite{daganzo84b} (Figure \ref{Strip_Strategy}), while it is claimed that the new strategy is more suitable for computer-based calculations and provides more flexibility with respect to the shape of the region. Numerical examples were conducted to demonstrate how this strategy can be applied to a VRP in a region of an irregular shape. 

This idea has been extended to the vehicle routing problem (VRP) where a fleet of vehicles (instead of a single salesman) serves region $\mathbf{T}$. Each customer has a random demand with the expectation of $v$. All vehicles shall start from and return to a depot $r$ distance away
from the center of $\mathbf{T}$, and each of them can serve $V$ unit demand. With a similar asymptotic analysis, it has been shown that the optimal value of the total expected all-vehicle travel distance
can be approximated by: 
\begin{equation}
2r(v/V)N+k_{VRP}\sqrt{TN},
\end{equation}
where $k_{VRP}$ is again a constant scalar determined by the metric. This approximated optimal distance can be achieved by using a ring and sweep strategy as shown in Figure \ref{Ring_Sweep_Strategy}. Traditional VRPs have been adapted and extended in several dimensions in recent two decades. With the rapid growth of VRP problem sizes, some studies restrict the VPR partition scheme to certain specific patterns (e.g., ring sectors) so as to cope with the challenges from increasing computational complexity. For these new problems, multiple distance metrics including the ring-radial metric have been explored. Several new factors that affect optimal routing design, such as service frequency and time windows, have been incorporated in systematic decisions.

In terms of traditional models, \citet{torres1999} consider a garbage pickup problem that requires all links in a transportation network to be visited which can be viewed as a generalization of the CPP by using multiple vehicles rather than a single postman for the delivery (or pickup). \citet{torres1999} proposed a CA model to create the routes for the garbage pickup problem in an urban region of the square grid metric. Built on grid metric approximation proposed by \cite{daganzo84b}, their CA approach includes a sequence of steps, i.e., identifying equi-travel time contours, determining zone sizes, laying out zone geometry, and finally using an iterative procedure to design the optimal routes.

\subsection{Distribution Studies}
\label{sec:distribution_studies}

In \citet{langevin96}, freight distribution systems are divided into six classes in terms of the general distribution structure and existence of transshipment: (1) one-to-many distribution without transshipment, (2) many-to-one distribution without transshipment, (3) many-to-many distribution without transshipment, (4) one-to-many distribution with transshipments, (5) many-to-many distribution with transshipments and lastly (6) integrated works. In the first three classes, the distribution is directly from the origin(s) to destination(s) without any transshipment in between whereas classes 4 and 5 include distribution centers (transshipment locations) over which the goods are distributed. In other words, the distribution is carried out in two stages for classes 4 and 5: The first stage is from origins to distribution centers and the second one is from distribution centers to the end users. The last class considers studies that integrate continuous and discrete models. The review paper does not have a many-to-one distribution with transshipments class since no such paper was written until 1996 (we have not found any recent studies under this setting either). Table \ref{tab:Distribution} provides a summary of problem types, cost components and main contributions of the distribution related studies in this review.

\setlength{\extrarowheight}{1em}
\begin{table}[!ht]
\footnotesize
 \caption{Distribution papers summary.}
 \vspace{0.25cm}
\begin{tabular}{c| L{4cm}| P{0.25cm} P{0.25cm} P{0.25cm} P{0.25cm} P{0.25cm} | P{7.5cm} }  
	& &\rotatebox[origin=l]{90}{ \hspace{-0.55cm} Routing} & \rotatebox[origin=l]{90}{ \hspace{-0.55cm} Vehicle related} & \rotatebox[origin=l]{90}{ \hspace{-0.55cm} Terminal}
   	& \rotatebox[origin=l]{90}{ \hspace{-0.55cm} Rent} & \rotatebox[origin=l]{90}{ \hspace{-0.55cm} Other} \\
   	&\multicolumn{1}{c|}{\textbf{Study}} & \multicolumn{5}{c|}{\textbf{Cost Structure}} & \textbf{Main Contribution(s)} \\
  	\toprule 
    	\multirow{10}{*}[-3em]{\begin{sideways}One-to-many\end{sideways}}  	
  	& Daganzo and Erera (1999) & X & X &   &   & X & \parbox{7.5cm}{\linespread{0}\selectfont Considers robust logistics systems with redundancy.} \\
  	& Erera and Daganzo (2003) & X &   &   &   &   & \parbox{7.5cm}{\linespread{0}\selectfont Extends Daganzo and Erera (1999) to dynamic VRP with uncertain demand.} \\
  	& Francis and Smilowitz (2006) & X &   &   &   &   & \parbox{7.5cm}{\linespread{0}\selectfont Introduces a CA model to PVRP with service choices.} \\
    & You et al. (2011)            & X &   &   &   & X & \parbox{7.5cm}{\linespread{0}\selectfont Integrates short-term distribution decisions with long-term inventory decisions in gas distribution systems.} \\
    & Agatz et al. (2011)          & X &   &   &   &   & \parbox{7.5cm}{\linespread{0}\selectfont Provides a fully automated approach that provides high-quality time slot schedules.} \\
	& Pang and Muyldermans (2012)  & X &   &   &   &   & \parbox{7.5cm}{\linespread{0}\selectfont Analyzes the impact of service postponement on routing costs.} \\    
    & Saberi and Verbas (2012)     & X &   &   & X & X & \parbox{7.5cm}{\linespread{0}\selectfont Facilitates strategic planning in a time-dependent environment while considering emission costs.} \\
	& Turkensteen and Klose (2012) & X &   &   &   &   & \parbox{7.5cm}{\linespread{0}\selectfont Derives logistics cost estimates from the dispersion of demand points for a one-to-many system.} \\    
    & Davis and Figliozzi (2013)   & X & X &   &   & X & \parbox{7.5cm}{\linespread{0}\selectfont Analyzes the competitiveness of electric vehicles considering energy consumption costs.} \\    
    & Huang et al. (2013)          & X &   &   &   &   & \parbox{7.5cm}{\linespread{0}\selectfont Introduces a hybrid strategy that improves upon the solution from CA with a tabu search in humanitarian setting.} \\
    & Lei et al. (2016)          & X& X& &  & X & \parbox{7.5cm}{\linespread{0}\selectfont Uses PVRP-SC concept in a bilevel game-theoretic model to improve the effectiveness of parking enforcement patrol.} \\
    \midrule    
    \multirow{5}{*}[-1em]{\begin{sideways}Many-to-many\end{sideways}}
	& Kawamura and Lu (2007)       & X &   & X & X &   & \parbox{7.5cm}{\linespread{0}\selectfont Analyzes the effect of consolidation in many-to-many systems.} \\
    & Smilowitz and Daganzo (2007) & X & X & X &   &   & \parbox{7.5cm}{\linespread{0}\selectfont Compares fully integrated networks and nonintegrated networks for expedited and deferred deliveries.} \\
    & Chen et al. (2012)           & X &   & X & X &   & \parbox{7.5cm}{\linespread{0}\selectfont Extends Kawamura and Lu (2007) to represent smaller businesses.} \\
    & Campbell (2013)              & X &   &   &   &   & \parbox{7.5cm}{\linespread{0}\selectfont Provides expressions for the optimal number of hubs, hub locations and the related cost in time definite setting.}\\
	& Lin et al. (2014)            & X &   & X & X &   & \parbox{7.5cm}{\linespread{0}\selectfont Extends the study of Chen et al. (2012) considering energy consumption and PM2.5 emissions as additional costs.} \\    
    \bottomrule
\end{tabular}
\label{tab:Distribution}
\end{table}

\citet{turkensteen2012demand} derive logistics cost estimates from the dispersion of demand points for a system which serves a set of target customers from a single central facility. Their study focuses on outbound costs and they assume that demands of the customers are almost identical and deterministic. Furthermore, the products are delivered with fixed routes at fixed intervals. The authors develop travel distance estimates and compare the estimates with distance approximations. They use CA to approximate travel distances for uniformly distributed demand points.

\citet{you2011optimal} consider a one-to-many distribution problem without transshipment in the scope of optimal distribution and inventory planning of industrial gases. The decision maker determines tank sizing (inventory profile) at each customer location in addition to delivery and routing related decisions with the aim of minimizing the total cost. The study assumes that a single type of industrial gas (single product) is distributed with trucks that have the same traveling speed. The authors propose a CA-based solution strategy that divides the problem into an upper and a lower level and solves them sequentially.

\citet{saberi2012continuous} provide a CA model to minimize emissions in the emission VRP (EVRP). EVRP is a variant of time-dependent VRP as emissions strongly depend on congestion and average travel speed which can change significantly between peak and off-peak periods. Their model aims to facilitate strategic planning of one-to-many distribution systems without transshipment in a time-dependent environment. The routing costs are mainly based on the VRP approximations in \citet{daganzo05}; however, related emission costs are also included at each stage. The problem is decomposed geographically into a set of sufficiently large sub-regions with nearly constant customer density. Instead of considering emission costs, \citet{davis2013methodology} use CA as a part of their integrated cost model to evaluate the competitiveness of electric delivery trucks. CA is used to approximate the average cost of serving routes in VRP setting. Assuming that the number of customers per route is balanced and each vehicle is capacitated, they use a refined version of VRP approximation of \citet{daganzo05}, which is introduced by \citet{figliozzi2008planning} with an additional term to modify the local (service region) tour distance. The authors state that the VRP cost term can be generalized to the cases with time window constraints by using the approximation of \citet{figliozzi2009planning} the average length of VRPs with time windows.

\citet{francis2006modeling} use a CA model to solve a periodic VRP with service choice (PVRP-SC), which is a VRP variant where the visit frequencies to the customers are considered as decisions of the model. This problem can be classified under one-to-many distribution without transshipment. The authors design a fast and valid approximation scheme for the strategic level PVRP-SC since the exact and heuristic solution methods are limited by the size of problem instance \citep{francis2006period}. Similar to VRP approximations in \citet{daganzo05}, the routing costs are divided into line haul and detour costs. The exact data for node locations and demand volumes are replaced with continuous density functions and a continuous model is developed for PVRP-SC. The concept of PVRP-SC is later used in a bilevel game-theoretic model to improve the effectiveness of parking enforcement patrol \citep{leiUNP}, where  individual parking drivers' endogenous payment decisions (which are based on knowledge of the patrol visit frequencies) directly affect the number of vehicles in parking violation as well as the patrol officers' work efficiency. The CA model is compared with a discrete mathematical programming formulation to reveal managerial insights and to show computational efficiency.  In a similar setting, \citet{pang2013vehicle} use CA to evaluate the value of postponing routing decisions on overall vehicle routing cost or distance when clients' requests accumulate over time and vehicles with fixed capacity visit clients regularly to cover the accumulated demand. Each client must be served with a single vehicle; i.e., demand splitting across vehicles is not allowed. They solve the underlying capacitated VRPs with a local search heuristic and CA. The CA model is also used in cases with different (asymmetric) demand rates. The results from local search and CA behave similarly and illustrate the benefit of postponement, yet the results of CA model overestimate the improvement compared to the local search as the asymmetric demand cases lead to inaccuracies in approximations.

%[*** reference: Lei, C., Zhang, Q. and Ouyang, Y. "Optimal patrol planning for urban parking enforcement considering drivers' payment behavior." Transportation Research Part B. Under review.]

\citet{agatz2011time} use CA to determine the set of delivery time slots to offer customers in different regions for an e-tailer in a one-to-many distribution system without transshipment. Similar to \citet{francis2006modeling}, the authors provide an exact integer programming alternative and test the validity and quality of their CA model on different instances. Different from \citet{francis2006modeling}, their CA model includes costs that stem from time slot configurations (such as costs within the same zip and costs across different zips in the same time slot). Their CA model is able to provide high quality solutions in significantly shorter times compared to those of integer programming.

In another one-to-many distribution without transshipment setting, \citet{huang2013continuous} use CA to solve the assessment routing problem in a disaster relief. CA becomes useful in a disaster relief setting due to lack of precise data, time pressure and necessity for easy implementation. The study considers two different routing policies based on the separation of the region of interest, extending \citet{newell1986design}. Different from already discussed studies, the authors propose a hybrid solution approach which improves the solutions of the CA model by using a restricted Tabu Search. Often CA methods are considered as complementary solution methods; this study provides a framework for integration of CA and exact and/or heuristic methods.

\citet{daganzo1999planning} study the effects of uncertainty on logistics system design. To design a robust logistics system, they determine the optimal level of redundancy, which is both sufficient and cost-effective, along with the operating strategy necessary to utilize this redundancy.  They consider two one-to-many problems, a static VRP (no transshipment) and a warehouse location-inventory-routing problem (with transshipment). Redundancy is introduced in the VRP by removing the single tour assumption made in the initial problem statement.  The delivery zones and secondary vehicle routes must be determined, which is done using CA.  They find that the total distance traveled is a function of primary and secondary line-haul distances as well as combined local delivery distance.  The primary vehicles follow their tours and return when they reach capacity.  The secondary vehicle tours serving the remaining customers can then be modeled deterministically. \citet{erera2003dynamic} further extend this model to dynamic vehicle routing with uncertain demand that is initially known only by its stochastic distribution and then realized over time. They suggest a threshold global sharing scheme, which requires the use of CA. In their model, the vehicles come from a central depot to serve the customers in the region.  The customer demand and locations are initially only given as a distribution.  As more information comes in, the vehicle routes are re-optimized in real-time in order to minimize cost.  With their scheme, the customers are served in three phases.  The first phase is a preplanned route.  When the vehicles reach capacity, they return to the depot, even if they have not completed their planned route.  When phase one is over, the remaining customers are assigned to the vehicles for phase two.  Phase three consists of any customers that are still remaining at the end of phase two.  The customer locations become known right before the vehicles leave the depot, and the demand size of each customer is only known when the vehicle reaches the customer. CA is used to estimate the total route length of all three phases.

In a many-to-many setting, \citet{Smilowitz2007} analyze integrated package distribution systems in a multiple mode, multiple service level package distribution networks. Assuming two transportation modes (air and ground) and two service levels (express and deferred), the authors study two different network configurations: Fully integrated and non-integrated. Demand data is considered deterministic and stationary; however, an extension for stochastic case is also provided. The authors develop a CA model which represents the solution as functions (densities) of location (variables) under certain assumptions. This study differs from the existing studies in the literature since it considers multiple transshipments and multi-stop peddling tours in many-to-many setting and it includes additional operating costs such as repositioning empty vehicles.

The studies of \citet{kawamura2007evaluation}, \citet{chen2012comparison} and \citet{lin2016sustainability} use CA when evaluating the effects of delivery consolidation. The studies consider many-to-many networks with multiple suppliers and customers. When the distribution is consolidated, the suppliers cooperate and send their products to consolidation centers through which these products are distributed whereas the distribution is carried out directly from suppliers to customers in the unconsolidated case. \citet{kawamura2007evaluation} analyze the effect of consolidation in urban areas. They assume that each customer is treated equally and the same number of customers are visited at each dispatch with identical vehicles of maximum capacity. Inventory costs are considered negligible. In addition, the consolidated setting assumes constant headways and coordination between the inbound and outbound shipments is not considered. With these assumptions continuous cost functions are developed by considering inbound, outbound and terminal costs, and the optimal costs for both settings are calculated for two alternative types of vehicles. \citet{kawamura2007evaluation}  conclude that consolidation is not significantly beneficial, particularly with the current regulations regarding vehicle size and weight. \citet{chen2012comparison} extend \citet{kawamura2007evaluation}  to represent smaller businesses. They estimate parameters using additional data source and base cost-effectiveness comparison to additional factors such as customer density, demand quantity and vehicle specifications. Furthermore, they consider coordinated delivery consolidation strategy with synchronized inbound and outbound trips in addition to direct delivery and uncoordinated consolidation strategies.  Their results show that consolidation can be cost-effective when the consolidation center is close to the suppliers, the underlying network consists of large number of customers and suppliers, and terminal operations are not expensive. \citet{lin2016sustainability} extend the study of \citet{chen2012comparison} considering energy consumption and PM2.5 emissions as additional costs. Their results suggest similar findings to \citet{chen2012comparison} for cost-effectiveness of consolidation strategies.

In another many-to-many distribution setting, \citet{campbell2013continuous} considers time definite freight transportation in a one dimensional service region. The objective of the study is to minimize the maximum travel distance in the hub network for each origin-destination pair in order to achieve a certain service level. The author compares two different settings to achieve a given level of service: Use of a sufficient number of optimally located hubs (optimal in terms of transportation costs) and use of (potentially) fewer hubs, which are sufficient to satisfy the service level, without considering the transportation cost. Assuming that the hubs are equally separated, the author provides a CA model to determine the locations of the hubs over the service region. In conclusion, studies of distribution related problems within the last two decades focus more on real life problems with operational constraints instead of generic problems such as the traveling salesman problem and the VRP.

Many of these models only yield continuous vehicle routing characteristics that are not suitable for implementation in practice. To enhance the practicality the CA modeling framework of stochastic VRP, \citet{ouyang2007design} formulates an automated approach to discretion such continuous solutions into discrete vehicle routing zones. The author achieves this by utilizing a disk model and weighted centroidal Voronoi tessellations.  Customer demand is modeled in a region with a ring-radial network.  The disk model is used to find an initial zone partition, including the approximate sizes and shapes of the partition areas. Then, the partition is further optimized by the weighted centroidal Voronoi tessellations, which equalizes the delivery load size in each zone.  To solve stochastic problems, the customer demand just needs to be re-defined after the exact discrete information becomes known.

\subsection{Transit Studies}
\label{sec:transit_studies}

%{\color{blue} Will be taken from Section 5 of Smilowitz Draft.}

CA has been used extensively in transit settings particularly for public transportation systems with one-to-many and many-to-one patterns; see \citet{langevin96} for a summary on the early papers. A critical issue in public transportation is to offer services that are relatively inexpensive and as competitive as automobiles in terms of transportation quality. For this purpose, many studies focus on designing efficient networks by using various types of service options and routing the vehicles in the existing systems more efficiently. Table \ref{tab:Transit} provides a summary of problem types, cost components and main contributions of the transit related studies considered in this review. In addition, CA is frequently used in traffic assignment problems to investigate traffic flow and route choice behaviors of travelers in a specific region (see \citet{jiang2011dynamic} for a recent review).

\setlength{\extrarowheight}{1em}
\begin{table}[!ht]
\footnotesize
\caption{Transit papers summary.}
\vspace{0.25cm}
\begin{tabular}{c| L{4cm} | P{0.25cm} P{0.25cm} P{0.25cm} P{0.25cm} P{0.25cm} P{0.25cm} | P{7cm} }  
	& & \rotatebox[origin=l]{90}{\hspace{-0.55cm} Operation} & \rotatebox[origin=l]{90}{\hspace{-0.55cm} Passenger} & \rotatebox[origin=l]{90}{\hspace{-0.55cm} Infrastructure}
   	& \rotatebox[origin=l]{90}{\hspace{-0.55cm} User} & \rotatebox[origin=l]{90}{\hspace{-0.55cm} Occupancy} & \rotatebox[origin=l]{90}{\hspace{-0.55cm} Disutility} \\
   	&\multicolumn{1}{c|}{\textbf{Study}} & \multicolumn{6}{c|}{\textbf{Cost Structure}} & \textbf{Main Contribution(s)} \\
  	\toprule
  	\multirow{3}{*}[-1em]{\begin{sideways}DRT\end{sideways}}
    & Diana et al. (2006)             &    \multicolumn{6}{c|}{\rule{2.4cm}{0.4pt}}  & \parbox{7cm}{\linespread{0}\selectfont Solves fleet sizing problem for DRT services with time windows with a probabilistic model.} \\
	& Chandra and Quadrifoglio (2013) & X &   &   &   &   & X & \parbox{7cm}{\linespread{0}\selectfont Estimates the optimal cycle length of DRT vehicles in feeder transit services.}\\ 
    & Li et al. (2016)                & X &   & X & X &   &   & \parbox{7.5cm}{\linespread{0}\selectfont Considers car sharing in electrical vehicles.} \\	
	\midrule	
	\multirow{2}{*}[-1em]{\begin{sideways}FRT\end{sideways}}
    & Munoz (2002)			& X &  &   &   &   &   & \parbox{7cm}{\linespread{0}\selectfont Considers driver contract and schedule.} \\
	& Ellegood et al. (2015)			& X &   &   &   &   &  & \parbox{7cm}{\linespread{0}\selectfont Analyzes the effect of mixed loading in school bus routing.} \\  	
	\midrule  	
  	\multirow{5}{*}[-3em]{\begin{sideways}Hybrid\end{sideways}}
  	& Aldaihani et al. (2004)         & X &  &  & X &   &   & \parbox{7cm}{\linespread{0}\selectfont Determines the number of zones allowing multiple on-demand vehicles at each zone and assuming the ride time on a fixed route is a function of zone count.} \\
    & Quadrifoglio et al. (2006)      &   \multicolumn{6}{c|}{\rule{2.4cm}{0.4pt}}   & \parbox{7cm}{\linespread{0}\selectfont Develops bounds on the maximum longitudinal velocity in MAST setting.} \\
    & Daganzo (2010)                  & X &   & X & X & X &   & \parbox{7cm}{\linespread{0}\selectfont Discusses network shapes and operating characteristics to make transit systems as attractive as using an automobile. Various different cost functions are used.} \\
    & Quadrifoglio and Li (2009)      &   &   &   & X &   &   & \parbox{7cm}{\linespread{0}\selectfont Derives expressions that estimate the critical demand densities to decide when to switch the operation from DRT to FRT and vice versa.} \\
    & Li and Quadrifoglio (2011)      &   &   &   & X &   &   & \parbox{7cm}{\linespread{0}\selectfont Expands Quadrifoglio and Li (2009) by considering two vehicle operation in each zone.} \\
    & Nourbakhsh and Ouyang (2012)    & X &   & X & X & X &   & \parbox{7cm}{\linespread{0}\selectfont Integrates the hybrid network structure of Daganzo (2010) with the flexible route idea of Quadrifoglio et al. (2006). Designs a flexible transit system for low demand areas.} \\
    & Ouyang et al. (2014)    & X &   &  & X & X &   & \parbox{7cm}{\linespread{0}\selectfont Extends Daganzo (2010) to design bus networks under spatially heterogeneous demand, allowing denser local routes for high demand areas with smaller spacings.} \\
    \bottomrule
\end{tabular}
\label{tab:Transit}
\end{table}

Transit related services are divided broadly into two service categories: fixed route transit (FRT) service and demand responsive transit (DRT) service. In FRT, a large-capacity vehicle follows a fixed route with a given schedule and passengers can get on the vehicle at predetermined stops. Due to its large capacity, the vehicle can consolidate many passengers, making FRT a cost-efficient option in public transportation. In DRT, the route is not fixed and relatively small vehicles provide a door-to-door transportation option for passengers (taxicabs, dial-a-ride). DRTs provide more flexibility to passengers yet they are more costly compared to FRTs.

\citet{munoz2002driver} studies an FRT system from a driver contract and schedule perspective. The objective cost function consists of vehicle operating cost and passenger delay cost. The author considers both contract scheduling work as well as the design of a timetable based upon pre-existing driver contracts. For the latter, the author uses CA in order to address the uncertainty of driver absences as well as to analyze the cost sensitivity of the final solution. The influence of driver's absence is analyzed by comparing the ideal model assuming that all drivers are present and the stochastic one accounting for driver absences. A numerical optimization method is then used to verify the optimal solution. %The study suggests a flexible contract where drivers still work forty hours a week, without splitting shifts, but do not have to work eight hours a day.

\citet{Diana06} consider a fleet sizing problem for DRT services with time windows and develop a probabilistic model which inputs the demand distribution over the service region along with the desired service quality. While determining time windows, the authors consider an allowable level of waiting time implied by the considered quality level. Their DRT setting allows consolidation; i.e., the vehicle is allowed to pick up additional passengers before dropping off the current one(s). Their model uses CA to estimate the distribution of leg lengths in vehicle routes where leg length specifies the distance between two consecutively visited points. In a similar DRT setting, \citet{li2016design} uses CA to design an electrical vehicle (EV) sharing system without consolidation. Different from regular car sharing, EV sharing system must consider nonlinear vehicle charging times. The authors decompose the corresponding region into smaller regions and solve the problem in sub-linear time with bisection algorithm. Decisions such as fleet sizes and sharing station locations are optimized jointly with EV sharing operations to ensure a reliable service level under stochastic and dynamic customer demand.

\citet{chandra2013model} develop a two-step queuing model which estimates the optimal cycle length of DRT vehicles in feeder transit services. Demand responsive feeder transit services connect passengers to major networks with on-demand vehicles that generally pick up and drop off multiple passengers in a single trip. Serving many passengers in each tour is beneficial as it increases ride sharing; however, it increases riding times for the passengers. Having too few passengers in each tour may not be practical especially when the demand of the service region is relatively high. The authors construct a utility based CA model by using inequalities that relate the number of served passengers to the cycle length. Similarly, \citet{ellegood2015continuous} uses CA to analyze the effect of mixed loading in school bus routing. School bus routing can be considered as an FRT and mixed loading indicates that a single vehicle can carry students from different schools at the same time. \citet{Aldaihani04}, \citet{Quadrifoglio06} and \citet{Quadrifoglio09} use CA to design hybrid networks which combine DRT and FRT. \citet{Aldaihani04} develop a model which divides the region into smaller zones such that across-zone transits are handled with fixed route vehicles whereas each zone is served by on-demand vehicles. Their model allows a maximum of three transfers for each passenger, assumes that each zone is covered by the same number of on-demand vehicles and prevents ride sharing on on-demand vehicles. The cost function tries to balance the cost to passenger (travel time) and cost of on-demand and fixed vehicles. Continuous approximations of travel times and passenger waiting times are used to get a closed form representation of the cost function which is proven to be convex for certain parameter values. Furthermore, their sensitivity analysis indicates that demand density affects the optimal number of zones, but not the number of fixed busses per route. \citet{Quadrifoglio06} study mobility allowance shuttle transit (MAST) in a rectangular region. In MAST, the fixed route vehicle is allowed to make deviations from the course to pick-up/drop-off passengers while maintaining an acceptable velocity level towards its destination. Their study focuses on developing bounds on the maximum longitudinal velocity of the MAST vehicle. Similar to Daganzo's strip strategy \citep{daganzo84b}, the authors find a sufficient condition on the longitudinal distance between consecutive demand points for which no-backtracking is optimal. Assuming no-backtracking and using CA, the authors find a closed form formula for the lower bound on maximum velocity as a function of region length, region width, demand density, service time (pick-up / drop-off) and average vehicle speed. In addition, they develop two upper bounds on optimal value of maximum velocity. A different hybrid is proposed in the study of \citet{Quadrifoglio09} for feeder lines that are used to connect service areas to a major transit network. The feeder lines work in FRT setting when the demand is high whereas they switch to DRT setting as the demand decreases. In this study, the authors identify the conditions that trigger the switch in a rectangular region. CA is used to formulate the utility functions of DRT and FRT with one-vehicle and two-vehicle settings, and switching conditions (critical demand densities) are identified based on various factors such as region geometry, vehicle speed and travel time. \citet{li20112} expand this work by considering a two-vehicle operation in each zone, with each vehicle serving half of the zone.  This model is better for a large service area with a high level of demand.

\citet{daganzo2010structure} focuses on designing network shapes and operating characteristics to make transit systems as attractive as using an automobile. Daganzo works on idealized square regions with uniform demand and provides network structures which combine grid and hub-and-spoke concepts. A hybrid of both networks can be formed by having a central square with grid structure in the center and branches from central square to cover periphery. Daganzo's cost model uses CA to approximate agency metrics such as infrastructure costs and vehicle travel costs; and user metrics such as expected number of transfers, waiting times, travel distance and time. The ideas are further extended in \citet{ouyang2014continuum} to design bus networks under spatially heterogeneous demand, where denser local routes are allowed for high demand areas with distinctly smaller spacings on a power-of-two scale. In a related study, \citet{nourbakhsh2012structured} integrate hybrid network structure of \citet{daganzo2010structure} with flexible route idea of \citet{Quadrifoglio06} to design a flexible transit system for low demand areas. The continuous formulation allows them to derive analytical conditions in which flexible-route transit provides the most economic passenger transportation service. In general, the applications of transit problems from public transportation and electrical vehicle perspective have been increased as alternatives to fuel consuming private cars become more desirable.

\section{Integrated Supply Chain and Logistics Studies}
\label{Sec_Hybrid_Model}

%{\color{blue} Will be taken from Section 3.4 from Smilowitz Draft and Section 4 of Ouyang Draft. Please see Table \ref{tab: Sample} for a sample table.(Should this section have subsections?)}

Integrated supply chain and logistics models have been widely taken into consideration recently. Whether it is the integration of location and inventory problems or the integration of location and transportation problems, the availability of better approximation tools to deal with the complexity can be identified as the main reason for this trend. See Table \ref{tab:Integrated} for the summary of integrated supply chain papers.

\setlength{\extrarowheight}{1em}
\begin{table}[!ht]
\caption{Integrated supply chain papers summary.}
\vspace{0.25cm}
  \centering
  \footnotesize
    \begin{tabular}{L{4cm} |P{0.25cm} P{0.25cm} P{0.25cm} P{0.25cm} P{0.25cm}| P{7cm}}
          & \begin{sideways} \hspace{-0.5cm} Facility opening \end{sideways} & \multicolumn{1}{c}{\begin{sideways} \hspace{-0.5cm} Facility operating \end{sideways}} & \begin{sideways} \hspace{-0.5cm} Transportation \end{sideways} & \multicolumn{1}{c}{\begin{sideways} \hspace{-0.5cm} Inventory \end{sideways}} & \begin{sideways} \hspace{-0.5cm} Other\end{sideways} &  \\
    \multicolumn{1}{c|}{\textbf{Study}}& \multicolumn{5}{c|}{\textbf{Cost Structure}} & \multicolumn{1}{c}{\textbf{Main Contribution(s)}} \\
    \toprule
    Erlebacher and Meller (2000)	& X     &       & X     & \multicolumn{1}{c}{X} &       & \parbox{7cm}{\linespread{0}\selectfont Considers location-inventory in a two-level distribution system.}\\
    Wang et al. (2004) 				& X     &       &       &       & X     & \parbox{7cm}{\linespread{0}\selectfont Studies the economics of park and ride facilities.} \\
    Shen and Qi (2007) 				&       &       & X     & \multicolumn{1}{c}{X} &       & \parbox{7cm}{\linespread{0}\selectfont Calculates the routing cost using CA.} \\
    Naseraldin and Herer (2008)	    &      & X     & X     & X      &      & \parbox{7cm}{\linespread{0}\selectfont  Uses CA approach for integrated supply chain design with first order inventory sharing.} \\
    Naseraldin and Herer (2011)		&      & X      & X     & X      & X     & \parbox{7cm}{\linespread{0}\selectfont Integrates location and inventory decisions allowing transshipments of inventory between facilities.} \\ 
    Mak and Shen (2012)				& X      & X       &      & X       &      & \parbox{7cm}{\linespread{0}\selectfont Considers dynamic sourcing and inventory sharing in supply chain design.} \\
    Tsao and Lu (2012) 				&       & \multicolumn{1}{c}{X} & X     & \multicolumn{1}{c}{X} &       & \parbox{7cm}{\linespread{0}\selectfont Integrates facility location and inventory allocation problem with transportation cost discounts.} \\
    Tsao et al. (2012) 				& X     &       & X     & \multicolumn{1}{c}{X} &       & \parbox{7cm}{\linespread{0}\selectfont Integrates a facility location-allocation with an inventory management problem.} \\
    Tsao (2013) 					& X     & \multicolumn{1}{c}{X} & X     & \multicolumn{1}{c}{X} &       & \parbox{7cm}{\linespread{0}\selectfont Studies distribution center network under trade credits.} \\
    Saberi and Mahmassani (2013)	&       &       & X     &       &       & \parbox{7cm}{\linespread{0}\selectfont Introduces an application of CA in airline hub location problem.} \\
    Pulido et al. (2015)			&       & \multicolumn{1}{c}{X} & X     & \multicolumn{1}{c}{X} &       & \parbox{7cm}{\linespread{0}\selectfont Considers the spatial and time dimensions in a unique logistics problem.} \\
    Tsao (2016)						&       &       & X     & \multicolumn{1}{c}{X} & X     & \parbox{7cm}{\linespread{0}\selectfont Discusses supply chain network design problems for deteriorating items with trade credits.} \\
    Tsao et al. (2016)				& X     &       & X     &       & X     & \parbox{7cm}{\linespread{0}\selectfont Considers multi-item distribution network design problems of multi-echelon supply chains under volume (weight) discounts on transportation costs.} \\
     Lim et al. (2016)				& X     &       & X      & X      &      & \parbox{7cm}{\linespread{0}\selectfont Considers the supply chain's agility.} \\
	\bottomrule
    \end{tabular}%
    \label{tab:Integrated}
\end{table}%

A pioneer in analytical study of integrated supply chain design is done by \cite{naseraldin2008integrating}. They model a continuous location-inventory problem on a line segment and analytically quantify the benefits of integrated modeling. Although this model considers only first-order inventory sharing, \cite{naseraldin2011location} extend the problem to allow transshipments of inventory between facilities. They point out several interesting results, e.g., that it can be optimal to increase the density of facilities when unit holding cost increases. Their results highlight the importance of modeling other forms of inventory sharing than the classical risk pooling. 

\cite{mak2012risk} follow the idea of \cite{lim2013} and further extend the CA model to incorporate inventory considerations into strategic facility location and network design decisions. Response to demand uncertainty is allowed via inventory sharing arrangements, over both the temporal dimension (e.g., through uncertainty pooling) and the spatial dimension (e.g., through responsive transshipments). It is shown how inventory management strategies and location decisions are jointly optimized to affect the facility number and spatial distribution.

\citet{erlebacher2000interaction} study a location-inventory model for designing a two-level distribution system. They assume the location of serving customers to be continuously distributed. More recently, \citet{tsao2012supply} study an integrated facility location and inventory allocation problem with transportation cost discounts. They use CA to reduce complexity of distance calculations and provide a high-level solution for the supply chain management (SCM) problem. This is the first study based on CA for a supply chain network design problem with transportation cost discounts and inventory decisions techniques. Similarly, \citet{tsao2012continuous} study an integrated facility location and inventory allocation problem for designing a distribution network with multiple distribution centers and retailers. CA model is used to represent this complex network. In this study, they not only locate the regional distribution centers and assign retail stores to them, but also determine the inventory policy at each location. In a later study, \citet{tsao2016multi} use CA to compare the effect of single cluster replenishment policy to joint cluster policy on supply chain network design. They assume that the locations of distribution centers and retailers are given and the transportation cost is affected by volume discounts. Their proposed solution defines the input data as continuous functions similar to \citet{tsao2012supply} and \citet{tsao2012continuous}. \cite{lim2016agility} study the design of a supply chain distribution network under demand uncertainty. Unlike traditional facility location models for supply chain design, they consider the supply chain's agility (i.e. ability to quickly respond to unexpected fluctuations in customer demand). They use CA approach to model customer demand such that demand points are evenly spread
on a two-dimensional infinite homogeneous plane. They find that it is optimal to decrease the density of distribution centers when transportation cost increases and to increase the density of distribution centers when the shortage penalty cost increases. Their findings suggest that proximity-based facility location models are not sufficient for designing modern responsive supply chains, and a new class of models is needed.

Similar to the location-allocation problem proposed by \citet{Murat10} but with trade credits, \citet{tsao2013distribution} considers a distribution center network design problem. Trade credit means that a seller gives a buyer extra time to pay or a discount for early payment as a credit. Tsao uses CA to model the problem and provides the location of the distribution centers (DCs), the allocation of retail stores to DCs and the joint replenishment cycle time at DCs. Adding the item perishability to the problem, \citet{tsao2016designing} uses CA to study an integrated facility location, inventory allocation and preservation effort problems by considering perishable items and trade credits. This study uses an approximation technique \citet{tsao2012continuous} to divide the network into smaller regions over which the discrete variable can be modeled using the slow varying functions.

Integrating pricing decision with facility location setting, \citet{wang2004locating} consider the optimal location for a park-and-ride facility and pricing of a park-and-ride facilities within a linear continuum region. As an application of CA in airline service network planning, \citet{saberi2013modeling} apply CA to the airline hub location and optimal market problem, more specifically in restricted and unrestricted single-hub systems. Unlike the airline service network planning literature, they model the optimal market for an airline in a competitive structure with multiple airlines and already located hubs. Combining routing and facility location problems, \citet{Shen07} use CA to estimate the transportation cost from facilities to customers in their strategic level facility location model. Even though this study is not considered as a facility location study, it uses the results derived from another CA model while solving a facility location problem. Considering time windows constraints, \citet{pulido2015continuous} introduce a CA-based methodology to study a combined location-routing problem. The methodology uses CA to determine the optimal warehouse density for the served region, assuming its customer density is known. The proposed methodology is similar to \citet{agatz2011time}, but no time slots are offered. In conclusion, applications of CA to more integrated settings have increased recently, due to the fact that the CA approach provides an easy-to-implement platform to deal with such complex problem settings. 

\section{Concluding Remarks and Future Directions}\label{Sec_Conclusions}

%OLD Conclusion: In this review, we survey the CA models that have been developed for facility location problems as well as transit and distribution problems. We aimed to classify the problems and to highlight the main contributions of each work along with how they used CA approach in their modeling and analysis. Moreover, we provide a guideline for the logistics researcher by addressing the fruitful topics and introducing the gaps need to be filled.
%%\noindent {\color{red} Basdere Note: Maybe change this section to `Future of Continuous Approximation'. We can start with concluding remarks and indicate that we are going to summarize current status of CA and then provide potential gaps for the future use of CA. Something like this...}

In this review, we survey the continuum approximation models that have been developed for facility location problems, distribution and transit problems, and integrated supply chains. We classify the problems and highlight the main contributions of each work along with a discussion on the use of CA in their modeling and analysis. In our concluding remarks, we provide a summary of the current status of CA and a guideline for logistics researchers by identifying the fruitful topics and introducing the gaps that need to be filled.

\subsection{Current Trends and Addressed Gaps}

Careful observation of the CA literature in the past two decades reveals two clear trends in CA applications.  First, the routing related problems were in the center of CA studies in the early 80s and 90s \citep{langevin96}. However, this trend seems to be shifting towards facility location related studies in the last two decades. Compared to the other streams, more researchers are focusing on the CA applications on facility location problems and integrated supply chains. Second, the research on CA is used for a wider range of applications (e.g. home grocery delivery, bus routing, electric delivery trucks, routing in disaster relief, traffic sensor deployment for surveillance, airline hub location and so forth) compared to the last decade. 

In their literature review, \cite{langevin96} discuss four main gaps in the CA literature which are investigated by various studies during the last decade: First, there were few reported applications of CA. CA approaches were used at early stages of the analysis and replaced by more detailed (discrete) models in later stages, which were reported in the literature instead of CA approaches. Along with many other application studies, the studies of \cite{Diana06}, \cite{Quadrifoglio09}, \cite{agatz2011time}, and \citet{you2011optimal} fill this gap with their reported applications on large scale real-life applications. Second, the cost structure of early studies was mostly limited in the sense that they only consider the trade-off between inventory and transportation costs. As the studies focus more on a wide range of real-life applications, various cost structures such as emission cost \citep{saberi2012continuous, lin2016sustainability}, energy cost \citep{davis2013methodology} and repositioning cost \citep{Smilowitz2007} are discussed. Third, the vast majority of the studies ignored time window or duration constraints which are essential components for many logistics systems. The studies of \cite{agatz2011time} and \cite{campbell2013continuous} fill this gap as they incorporate time notion under various settings. Last, certain logistics systems (many-to-many distribution with transshipment and peddling) had not been studied. The integrated package distribution study by \cite{Smilowitz2007} and the consolidation studies by \citet{kawamura2007evaluation}, \citet{chen2012comparison} and \citet{lin2016sustainability} are examples of many-to-many distribution with transshipment systems.

\subsection{Research Gaps}

The rest of this section discusses potential research topics and gaps in the existing literature. These gaps can be clustered under four main topics: 

\begin{enumerate}[label=(\roman*)]
\item \textbf{Using CA to support exact or heuristic solutions procedures:} In most studies, CA approach is used as an alternative to hard-to-solve discrete problems. These studies chiefly discuss the difficulty of the exact model and elaborate on how to use a CA model to yield feasible solutions within reasonable time. Nonetheless, CA is not only an effective tool to generate fast solutions, but also it can be used to generate upper and lower bounds for the original problem as in \cite{ouyang06}. Such bounds enable decision makers evaluate the quality of the solutions and have potential to improve the speed and quality of the exact models. In other words, rather than being an alternative approach, CA can be integrated to exact solution methods.
\item \textbf{Using CA at the operational level:} CA approach is chiefly applied to make strategic level decisions and gain insights about the big picture. However, as seen in \citet{agatz2011time}, CA can be used in real-time tactical and even in operational problems to increase the responsiveness in a dynamic decision making environment. Most of the dynamic decisions have to be made in an order of minutes or seconds and solely using exact methods are likely to fail in these cases. Therefore, use of CA in dynamic decision making environment can be a potential future research topic.
\item \textbf{Designing better quality measures to evaluate CA solutions:} Apart from the performance or optimality gaps, we need better measures to quantify the solution quality of CA. Solutions of many CA models provide significant improvements with simple modifications. Further improvements can be achieved by solving exact models which can be extremely difficult \citep{daganzo05}. In addition to their difficulty, the solutions returned by exact models may be difficult to interpret and implement. From this perspective, CA models have a critical benefit as they generally lead to simpler solutions that are easy to understand, explain and implement. In other words, easiness in the implementation stage has a value in real life. \citet{huang2013continuous} is an example of such studies where the authors use Tabu search on the result of CA model to reduce the performance gap of the solution.  For this reason, defining measures which are capable of quantifying the trade-off between exact solutions and CA solutions in terms of both solution quality and ability to implement is necessary for a fair evaluation of both approaches.
\item \textbf{Extending the use of CA in integrated problems and other fields:} Integrated models have always been of interest for supply chain and logistics researchers because of their comprehensiveness, but they are not always tractable due to their complexity in both modeling and solution approaches. However, the CA approach may potentially bridge this gap. There are several studies addressing this gap with the integration of different aspects of supply chain and logistics modeling especially with facility location decisions \citep{tsao2012supply, tsao2012continuous, tsao2016multi, tsao2016designing, lim2016agility}. Furthermore, the application domain of the CA can be extended to other domains such as wireless sensor networks which satisfy most assumptions of CA.
\end{enumerate}

\section*{Acknowledgement}
We woud like to thank Chrisitine Rhoades helped with collecting and summarizing some of the references. 
\newpage

\section*{References}

\bibliography{mybibfile}

\end{document}